\documentclass[12pt]{amsart}

\usepackage[matrix,arrow,graph]{xy}
\swapnumbers
\numberwithin{equation}{subsection}

\theoremstyle{plain}
\newtheorem{thm}[equation]{Theorem}
\newtheorem{lemma}[equation]{Lemma}
\newtheorem{prop}[equation]{Proposition}
\newtheorem{cor}[equation]{Corollary}
\theoremstyle{definition}

\newtheorem{defn}[equation]{Definition}

\newtheorem{blank}[equation]{}

\newtheorem*{rmk}{Remark}

\newcommand{\la}{\langle}
\newcommand{\ra}{\rangle}

\newcommand{\wt}{\widetilde}
\newcommand{\wh}{\widehat}
\newcommand{\ol}{\overline}
\newcommand{\half}{\frac{1}{2}}

\def\id{{\rm id}}

\def\Hom{\mathop{\rm Hom}\nolimits}
\def\hom{\mathop{\rm hom}\nolimits}

\def\End{\mathop{\rm End}\nolimits}
\def\Spec{\mathop{\rm Spec}}

\newcommand{\C}{{\mathbb C}}

\newcommand{\Z}{{\mathbb Z}}
\newcommand{\R}{{\mathbb R}}
\newcommand{\kk}{{\R}}

\newcommand{\PP}{{\mathbb P}}

\newcommand{\cal}{\mathcal}

\newcommand{\cS}{{\cal S}}

\newcommand{\cI}{{\cal I}}
\newcommand{\cL}{{\cal L}}
\newcommand{\cP}{{\cal P}}
\newcommand{\cA}{{\cal A}}
\newcommand{\cR}{{\cal R}}
\newcommand{\cD}{{\cal D}}

\newcommand{\cO}{{\cal O}}

\newcommand{\sig}{{\sigma}}
\newcommand{\Sig}{{\Sigma}}
\newcommand{\udot}{{\scriptscriptstyle \bullet}}

\newcommand{\Ext}{\mathop{\rm Ext}\nolimits}

\newcommand{\ext}{\mathop{\rm ext}\nolimits}
\newcommand{\mof}{\text{\rm -mod}}
\newcommand{\lf}{{\rm lf}}
\newcommand{\Mof}{\text{\rm -Mod}}

\newcommand{\cRHom}{\mathop{R{\cal H}{\it om}}\nolimits}

\newcommand{\Lie}{\mathop{\rm Lie}}

\newcommand{\bD}{{\mathbf D}}

\newcommand{\bP}{{\mathbf P}}

\newcommand{\Lu}{L^\udot}
\newcommand{\Nu}{N^\udot}
\newcommand{\Iu}{I^\udot}

\newcommand{\LDelt}{{{\Lambda}_\Delta}}
\newcommand{\LSig}{{{\Lambda}_\Sig}}

\newcommand{\LDmod}{{\LDelt}\mof}
\newcommand{\LSmod}{{\LSig}\mof}
\newcommand{\Lsmod}{{\Lambda_{[\sig]}\mof}}
\newcommand{\Lsvmod}{{\Lambda_{[\sig^\vee]}\mof}}
\newcommand{\real}{\mathop{\rm real}\nolimits}
\newcommand{\st}{\natural}
\newcommand{\vspan}{\mathop{\rm span}}

\newcommand{\codim}{\mathop{\rm codim}}

\begin{document}
\title{Koszul duality for toric varieties}

\author{Tom Braden}
\address{Dept.\ of Mathematics and Statistics\\
         University of Massachusetts, Amherst}
\email{braden@math.umass.edu}
\date{\today}

\begin{abstract}
We show that certain categories of perverse sheaves on affine
toric varieties $X_\sig$ and $X_{\sig^\vee}$ defined by dual cones
are Koszul dual in the sense of Beilinson, Ginzburg and Soergel \cite{BGS}.
The functor expressing this duality is constructed explicitly
by using a combinatorial model for mixed sheaves on toric varieties.

\end{abstract}

\thanks{This work was supported in part by NSF grant DMS-0201823 }

\maketitle

\hyphenation{Mac-Pher-son}
\hyphenation{canon-ical}

\section{Introduction}

In \cite{BGS} Beilinson, Ginzburg, and Soergel described a remarkable 
duality on the category of highest-weight modules for a semisimple
complex Lie algebra $\mathfrak g$.  Expressed algebraically, it says
that the ring governing this abelian category, endowed with a certain
grading, is Koszul self-dual.  Koszul duality is a relation between 
graded rings which generalizes the duality
between symmetric and exterior algebras.

By localizing highest weight modules, this can be expressed as
a duality on the category of Bruhat-constructible (mixed)
perverse sheaves on the flag variety associated to $\mathfrak g$.
However, the methods of \cite{BGS} do not completely explain how Koszul 
duality might arise from geometry.  In fact, no combination 
of standard functors such as $Rj^*$ and $j^!$ can have the right 
effect on the mixed structure.
    
In this paper we construct an analogous Koszul duality between 
categories of perverse sheaves on a pair of affine toric varieties defined 
by dual cones.  Our proof is ultimately combinatorial, but it is more
``local'' than the proof of \cite{BGS}, and we hope that our methods
will eventually lead to a more geometric understanding of Koszul duality.
 
Although the singularities of toric varieties are in many ways 
simpler than the ones appearing in the Bruhat stratification,  
the two situations satisfy very similar combinatorics.  
For instance, Stanley \cite{St} proved a
convolution identity for $g$-polynomials of convex polytopes 
(which measure local intersection cohomology of toric varieties) 
in which both a polytope $P$ and its polar $P^\vee$ appeared, along
with their faces.  This identity is exactly analogous to an identity for 
Kazhdan-Lusztig polynomials from \cite{KL}.  As was noted in \cite{BGS},
Koszul duality for highest weight modules can be seen as a categorical 
lifting of the Kazhdan-Lusztig identity.  This 
was our original motivation to conjecture that there would be a similar 
duality for toric varieties.

The appearance of dual cones and the strange interaction with 
(mixed) Hodge structures suggests that there should be some
relation with mirror symmetry.  We do not yet understand what 
such a relation might be, but we remark that Stanley's 
identity was used in \cite{BatBor} in the proof of the equality of
stringy Hodge numbers for mirror dual toric hypersurfaces.

\subsection{Main result, first version}\label{first main section}

Let $V$ be an $n$-dimensional real vector space, spanned by a 
lattice $V_\Z \subset V$.  Let $\sig\subset V$ be an $n$-dimensional pointed
polyhedral cone, rational with respect to $V_\Z$.  We denote by 
$[\sig]$ the poset (fan) of faces of $\sig$, with order relation $\prec$
given by inclusion of faces.  

The dual cone $\sig^\vee \subset V^*$ is defined by
\[\sig^\vee=\{\xi \in V^*\mid \xi(v) \ge 0 \text{ for all } v\in \sigma\}.\]
There is an order reversing bijection $\tau \mapsto \tau^\bot$ between
faces of $\sig$ and faces of $\sig^\vee$.

Let $X_{[\sig]}$ denote the affine toric variety associated to 
the fan $[\sig]$.  The action of the complex
torus $T = \Hom(V^*_\Z, \C^*)$ on $X_{[\sig]}$ has  
one $T$-orbit $O_\tau$ for each face $\tau$ of
$\sig$, so that $O_\rho \subset \ol{O_\tau}$ if and only if $\tau \prec \rho$.

In \S\ref{intro to pcs} below we define a full abelian 
subcategory $\cP_\Phi = \cP_\Phi(X_{[\sig]})$ of  
the category of perverse sheaves on $X_{[\sig]}$.  It depends 
on certain auxiliary data $\Phi$ attached to $[\sig]$, which we call
a ``combinatorial completion'', or ``completion'' for short.  
The simple objects of $\cP_\Phi$ are $\{L_\tau\}_{\tau\in[\sig]}$, where 
$L_\tau = IC^\udot(\ol{O_\tau}; \R)$ is
the intersection cohomology sheaf (with constant coefficients) 
supported on $\ol{O_\tau}$.  The completion $\Phi$ gives restrictions
on the allowed extensions between objects of $\cP_\Phi$; one consequence
is that $\cP_\Phi$ has enough projectives and enough injectives.

A completion $\Phi$ of $[\sig]$ induces a dual completion $\Phi^\vee$
of $[\sig^\vee]$, so we get a category of perverse sheaves
$\cP_{\Phi^\vee} = \cP_{\Phi^\vee}(X_{[\sig^\vee]})$ with simple
objects $\{L_{\tau^\bot}\}_{\tau\in[\sig]}$.
Let $I_\tau \in \cP_\Phi$ and $I_{\tau^\bot}\in \cP_{\Phi^\vee}$ 
denote the injective hulls of
$L_\tau$, $L_{\tau^\bot}$, so that
$I = \bigoplus_{\tau\prec\sig} I_\tau$ and 
$I^\vee = \bigoplus_{\tau\prec\sig} I_{\tau^\bot}$ are injective 
generators of $\cP_\Phi$ and $\cP_{\Phi^\vee}$, respectively.
Also set $L = \bigoplus L_\tau$, $L^\vee = \bigoplus L_{\tau^\bot}$.

The first version of our main result is the following.  It is exactly
analogous to the main result of \cite{S}, where the 
categories $\cP_\Phi$ and
$\cP_{\Phi^\vee}$ take the place of perverse sheaves on the flag variety.

\begin{thm}\label{first main}
There is a canonical isomorphism of rings
\[\End_{\cP_{\Phi^\vee}}(I^\vee)^{opp} 
\cong \Ext^\udot_{\cP_{\Phi}}(L, L).\]
The ring on the right, with its natural grading, is Koszul.
\end{thm}

Let us explain this result in more detail.
For a graded ring $A$, let  $A\mof$, $A\Mof$ denote the 
categories of finitely generated graded and ungraded left $A$-modules, 
respectively. 
Let $A$, $A^\vee$ be the endomorphism rings of 
$I$, $I^\vee$, 
respectively, given gradings via Theorem \ref{first main}. 
The functor $\Hom(-, I)$ gives an equivalence of 
categories $\cP_{\Phi}^{opp}\to A\Mof$, which sends
$L$ to $A_0 = A/A_{>0}$.

A graded ring $A = \bigoplus_{k\ge 0} A_k$ is called Koszul if
(1) $A_0$ is semisimple, and (2) $\Ext^i_{A\mof}(A_0, A_0\la -j\ra) = 0$
unless $i=j$, where $\la j\ra$ is the functor of shifting the grading
down by $j$ (note that this is the opposite convention to the one in 
\cite{BGS}).

The theorem can thus be interpreted as saying 
that $A^\vee$ is the opposite ring to $\Ext^\udot_{A\Mof}(A_0,A_0)$.
Rings related in this way are said to be Koszul dual to each other.


\subsection{Main result, second version}  Theorem \ref{first main}
can also be expressed by the existence of a duality functor between derived
categories of $A$-modules and $A^\vee$-modules.
Let $\wt L_\tau$ denote the simple object in $A\mof$
concentrated in degree $0$
whose ungraded version is $L_\tau$.  Let $\wt I_\tau$ be its injective hull;
it will be a graded version of $I_\tau$.  
Denote the simple and indecomposable injective objects of
$A^\vee\mof$ similarly.

\begin{thm} \label{Koszul functor} There is a contravariant 
equivalence of triangulated categories
\[\kappa = \kappa_\sig\colon D^b(A\mof)^{opp} \to 
D^b(A^\vee\mof),\]  
for which 
 \begin{enumerate}
\item[(a)]  
$\kappa(M[k]) = (\kappa M)[-k]$ and  $\kappa(M\langle k\rangle) = 
(\kappa M)[-k]\langle k \rangle$ for any $M \in  D^b(A\mof)$, and 
\item[(b)] 
$\kappa(\wt{L}_\tau) \cong \wt I_{\tau^\bot}$ and
$\kappa(\wt I_\tau) \cong \wt L_{\tau^\bot}$ for any $\tau\prec \sig$.
\end{enumerate}
\end{thm}

In \cite{BGS} such a functor $\kappa$ is constructed for
any pair $(A,A^\vee)$ of Koszul dual rings.  We take the opposite point of
view, however, first constructing $\kappa$ which satisfies
Theorem \ref{Koszul functor}, and then deducing 
Theorem \ref{first main} by the following argument.

Put $\wt L = \bigoplus \wt L_\tau \cong A_0$ and $\wt I^\vee = \bigoplus 
\wt I_{\tau^\bot}$.  
Then Theorem \ref{Koszul functor}(b) says 
that $\kappa(\wt L) = \wt I^\vee$.  If we define $\hom^\udot(M,N) = 
\bigoplus_{i,j\in \Z} \Hom(M, N[i]\la j \ra)$, then 
$\kappa$ induces a (non-grading preserving) ring isomorphism 
\[\hom^\udot_{D^b(A^\vee\mof)}(\wt I^\vee,\wt I^\vee)^{opp} \cong 
\hom^\udot_{D^b(A\mof)}(\wt L, \wt L).\]
Since higher $\Ext$'s between injective objects vanish, the left hand side 
is just $A^\vee$; the isomorphism of Theorem \ref{first main} follows.
Finally, for any $\tau,\rho \in [\sig]$ we have 
\begin{equation*}
\begin{split}
\Ext^i_{A\mof}(\wt L_\tau,\wt L_\rho\la -j\ra) & = 
\Hom_{D^b(A\mof)}(\wt L_\tau,\wt L_\rho\la -j\ra[i]) \\
  & =  \Hom_{D^b(A^\vee\mof)}(\wt I^\vee_\rho, \wt I^\vee_\tau[i-j]\la j\ra)\\
  & = \Ext^{i-j}_{A^\vee\mof}(\wt I^\vee_\rho, \wt I^\vee_\tau\la j\ra).
\end{split}
\end{equation*}
This vanishes if $i\ne j$, which proves that $A$ is Koszul.

\subsection{Combinatorial completions}\label{intro to pcs}
Our Koszul duality does not apply directly to 
the category of orbit-constructible perverse sheaves on toric varieties,
because there are not enough injectives; there
are local systems with arbitrarily long composition
series on every orbit, so the same is true of perverse sheaves.  
We solve this problem by restricting to a subcategory 
of perverse sheaves with prescribed monodromy, in the following way.  

Let $\Delta$ be a rational fan.  More generally, 
we can let $\Delta$ be the set of cones of a fan which 
lie outside a subfan;  we call such a set
a ``quasifan''.  A rational quasifan corresponds to a locally closed 
union of orbits in a toric variety, which we still denote $X_\Delta$.  

For any cone $\sig\in \Delta$, the $T$-orbit 
$O_\sig \subset X_\Delta$ has a unique smallest $T$-invariant 
neighborhood $U_\sig = \bigcup_{\tau\prec\sig} O_\tau$.  It is
isomorphic to $O_\sig \times X_{[\hat\sig]}$,
where $X_{[\hat\sig]}$ is the affine toric variety corresponding to 
the cone $\hat\sig$, which is equal to $\sig$ as a set, but 
considered as a cone in $V_\sig = \vspan(\sig)$.  This isomorphism
is not canonical, however.  Although there is a canonical projection
$U_\sig \to O_\sig$, there is some freedom in choosing 
a projection $\pi_\sig\colon U_\sig \to X_{[\hat\sig]}$. 

Choosing the projection $\pi_\sig$ is equivalent to fixing  
a subtorus of $T$ complementary to the stabilizer 
$T_\sig$ of any point in $O_\sig$.  A combinatorial completion 
should be thought of as 
the choice of such a complement for every $\sig$ in a compatible way.
To put this in terms of fan geometry, we use the 
identification $\Lie T\cong V\otimes_\R\C$, which takes 
$\Lie T_\sig$ to $V_\sig\otimes_\R\C$.  
\begin{defn} A combinatorial completion 
$\Phi$ of a quasifan $\Delta$ is a choice of a subspace 
$\Phi_\sig\subset V^*$ for every $\sig\in \Delta$ 
so that 
\begin{enumerate}
\item if $\tau \prec \sig$, then $\Phi_\tau \subset \Phi_\sig$, and
\item $V^* = \Phi_\sig\oplus V^\bot_\sig$ for every $\sig \in \Delta$,
where $V^\bot_\sig$ denotes the annihilator of $V_\sig$.
\end{enumerate}

If, in addition, $\Delta$ is rational, and
\begin{enumerate}
\item[(3)] $\Phi_\sig \cap V^*_\Z$ generates $\Phi_\sig$
 for every $\sig \in \Delta$
\end{enumerate}
(here $V^*_\Z \subset V^*$ is the lattice dual to $V_\Z$), then we say that 
$\Phi$ is rational.  
\end{defn}
The choice of $\Phi$ induces projection maps
$\pi_\sig\colon U_\sig \to X_{[\hat\sig]}$ for $\sig \in \Sig$; 
see \S\ref{Phi-stable}.
Properties (1)--(3) are actually slightly weaker than the geometric idea 
we started with; they only imply that the resulting  
map $U_\sig \to O_\sig\times X_{[\hat\sig]}$ is finite,
rather than an isomorphism.  However, this is sufficient for our purposes. 

Now we can define the category $\cP_\Phi(X_\Delta)$ referred to in 
\S\ref{first main section}.
\begin{defn} \label{defining stable perverse sheaves}
Let $\Phi$ be a rational combinatorial completion of a 
rational quasifan $\Delta$.  Define $\cP_\Phi(X_\Delta)$ to be the
full subcategory of $\cP(X_\Delta)$ consisting of objects 
$P$ so that for every $\sig\in \Delta$ there exists an isomorphism 
$P|_{U_\sig}\cong \pi_\sig^*P_\sig[\codim(\sig)]$, where
$P_\sig$ is a perverse sheaf on $X_{[\hat\sig]}$.  Such a $P$ will necessarily 
be constructible for the orbit stratification.
\end{defn}

\begin{rmk}  
Another, perhaps less artificial, way to
find a category of perverse sheaves with enough injectives
is to take limits of perverse sheaves whose simple constituents 
are $L_\tau$, $\tau\in [\sig]$.  The resulting objects
have unipotent monodromy on each orbit.  They are
not constructible in the classical sense,
as the stalks may be infinite dimensional, but a satisfactory
theory does exist.  The category on the Koszul dual side is then
(a mixed version of) the $T$-equivariant derived category.  
This is explained in \cite{BL}.
\end{rmk}

\subsection{Combinatorial mixed sheaves}
We construct our functor $\kappa$ by replacing the 
categories $D^b(A\mof)$ and $D^b(A^\vee\mof)$ by
equivalent categories $\bD_\Phi([\sig])$ and 
$\bD_{\Phi^\vee}([\sig^\vee])$.  These categories
provide a combinatorial model for mixed sheaves
on toric varieties, and should be of independent interest.
Two other flavors of combinatorial mixed sheaves appear
in \cite{BL}.

Given a fan or quasifan $\Delta$ with a completion $\Phi$ we define
triangulated categories $\bD(\Delta)$ and $\bD_\Phi(\Delta)$;
objects of these categories are complexes of sheaves
on the finite topological space $\Delta$. 
The definition of these categories 
is valid even for non-rational fans, when the toric variety does not exist.

We also define the following additional structures relating 
these categories:
\begin{itemize}
\item a ``forgetful functor'' $F_\Delta\colon \bD_\Phi(\Delta)\to \bD(\Delta)$,
\item pushforward and pullback functors $j^*$, $j^!$, $j_*$, $j_!$ 
defined for any inclusion $j$ of quasifans, and commuting with the
forgetful functors,
\item ``twist'' functors $\la k\ra$, $k\in \Z$, on both categories, commuting
with $F_\Delta$, and 
\item $t$-structures on both categories, with
corresponding abelian subcategories of perverse objects $\bP_\Phi(\Delta)$
and $\bP(\Delta)$.
\end{itemize}
$F_\Delta$ is compatible with these $t$-structures, in the sense that
the $t$-structure on $\bD_\Phi(\Delta)$ is pulled back from the 
$t$-structure on $\bD(\Delta)$.  In particular we have 
$\bP_\Phi(\Delta) = F_\Delta^{-1}(\bP(\Delta))$.  The twist functors $\la k\ra$ are also 
$t$-exact, and so they induce functors on $\bP_\Phi(\Delta)$ and
$\bP(\Delta)$.

Up to isomorphism and 
twists, there is a unique simple object $L^\udot_\tau$
in $\bP_\Phi(\Delta)$ which is supported on the closure of $\tau$. 
$F_\Delta$ induces a bijection between the isomorphism classes of 
simple objects of $\bP_\Phi(\Delta)$ and of $\bP(\Delta)$.

We also show that $F_\Delta$ embeds $\bP_\Phi(\Delta)$ as a full subcategory
of $\bP(\Delta)$. 
\begin{prop} \label{PC prop}  
$\bP_\Phi(\Delta)$ has enough injectives, and we have an 
 equivalence of categories 
$D^b(\bP_\Phi(\Delta)) \cong \bD_\Phi(\Delta)$.
\end{prop}

If the (quasi)fan $\Delta$ is rational and $X_\Delta$ is the corresponding toric variety, 
we define in \S5 a topological 
realization functor \[\real_\Delta\colon \bD(\Delta) \to D^b_u(X_\Delta),\]
where $D^b_u(X_\Delta)$
is the subcategory $D^b(X_\Delta)$ 
consisting of objects whose cohomology sheaves restricted to 
any orbit are local systems with unipotent monodromy.
We conjecture that this can be factored through a ``mixed''
realization functor sending $\bD(\Delta)$ to mixed Hodge modules
or mixed $l$-adic sheaves on $X_\Delta$.

\begin{thm} \label{realization theorem}
The realization functor $\real_\Delta$
satisfies:

\begin{enumerate}
\item For each $k\in \Z$, there is a natural isomorphism
\[\real_\Delta \circ \la k\ra \simeq \real_\Delta.\]
\item For any objects $S_1^\udot, S_2^\udot \in \bD(\Delta)$, 
the induced map
\[\bigoplus_{k\in \Z} \Hom_{\bD(\Delta)}(S_1^\udot, S_2^\udot\la k\ra)
\to \Hom_{ D^b_u(X_\Delta)}(\real_\Delta S_1^\udot, 
\real_\Delta S_2^\udot)\]
is an isomorphism.
\item If $j\colon \Sig\to \Delta$ is an inclusion of quasifans,
there exist natural isomorphisms
$j^* \circ \real_\Delta \simeq \real_\Sig \circ j^*$, 
$Rj_* \circ \real_\Sig \simeq \real_\Delta \circ j_*$, and similarly
for $j^!$, $j_!$.  Here we use the same letter $j$ to denote
the inclusion $X_\Sig\to X_\Delta$.
\item An object $S^\udot\in \bD(\Delta)$ is in $\bP(\Delta)$ if and 
only if $\real_\Delta S^\udot$ is a perverse sheaf.
\item $\real_\Delta|_{\bP(\Delta)}$ sends simples to simples.
\end{enumerate}
\end{thm}

If $\Delta$ is a rational quasifan with rational completion $\Phi$,
let $\real_{\Delta,\Phi}$ 
denote the composition \[
\bD_\Phi(\Delta)\to \bD(\Delta) \stackrel{\real_\Delta}{\longrightarrow}
D^b_u(X_\Delta);\]
it restricts to an exact functor $\bP_\Phi(\Delta) \to \cP(X_\Delta)$.

\begin{thm} \label{PC theorem}\indent\par 
\begin{enumerate}
\item If $S^\udot \in \bP(\Delta)$, then  $S^\udot$ is in 
$\bP_\Phi(\Delta)$ if and only if 
$\real_\Delta(S^\udot)$ is
in $\cP_\Phi(\Delta)$.
\item 
$\real_{\Delta,\Phi}$ takes injective objects
of $\bP_\Phi(\Delta)$ to injective objects in $\cP_\Phi(\Delta)$,
and all injectives in $\cP_\Phi(\Delta)$ can be 
obtained in this way.
\item For any two objects $S_1^\udot, 
S_2^\udot \in \bP_\Phi(\Delta)$, and any $i\ge 0$,
we have
\[\Ext^i_{\cP_\Phi}(\real_{\Delta,\Phi}S_1^\udot, \real_{\Delta,\Phi}
S_2^\udot) \cong
\bigoplus_{k\in\Z} \Ext^i_{\bP_\Phi(\Delta)}(S_1^\udot,S_2^\udot\la k \ra).\]
\end{enumerate}
\end{thm}

Let us explain the intuition behind the construction of 
$\bD(\Delta)$.  There is a  ``combinatorial
moment map'' $\mu\colon X_\Delta \to \Delta$ which sends all points
of $O_\sig \subset X_\Delta$ to $\sig$; put another way, it is the quotient
map $X_\Delta \to X_\Delta/T$.  Then  
$\Lambda = \bigoplus_j R^j\mu_*\R_{X_\Delta}$ is a sheaf of graded rings
on $\Delta$, whose stalk $\Lambda_\sig$ 
at a face $\sig$ is 
the cohomology ring $H^\udot(O_\sig;\R)$, which is an 
exterior algebra.

The cohomology of an object $S \in D^b(X_\Delta)$ 
is naturally a $\Lambda$-module.  More 
information is preserved if we take the derived push-forward;
$R\mu_*S$ becomes an object in the category of $dg$-modules over
$\Lambda$, considered as a sheaf of differential graded algebras 
with trivial differential.  This works because 
the torus $T$ is formal, meaning that the 
$dg$-algebra $\Omega^\udot(T)$ of differential forms is quasi-isomorphic
to its own cohomology.  A similar idea appears in Lunts \cite{L}.

Mixed complexes are ``bigraded'' objects: there is one shift operation
on degrees, and another on weights.  Accordingly, objects of
our category $\bD(\Delta)$ are complexes of graded $\Lambda$-modules 
rather than $dg$-modules.  This gives two shift operators, one
shifting the complex degree and another shifting the grading.
The $t$-exact twist functor $\la k \ra$ of Theorem \ref{realization theorem}
is a diagonal combination of both types of shift.

\begin{rmk}
A key reason we can describe weights so simply is that for each $i$ the 
cohomology $H^i(O_\tau;\R)$ is pure (of weight $2i$).
To make a similar construction on flag varieties, one wants to 
replace the exterior algebras $\Lambda_\sig$ by the stalk cohomology
groups of $j_{S*}\R_S$, where $S$ is a Bruhat cell and $j_S$ its inclusion
into the flag variety.  Unfortunately this cohomology is not pure, as
follows from an example of Boe (\cite{Boe}, see also \cite{SSV}).  
We hope that it will still be possible to understand 
Koszul duality on flag varieties using these 
techniques, perhaps by staying on the level of $dg$-algebras.
\end{rmk}

\subsection{Main result, third version}
If $\Delta = [\sig]$ for a rational fan $\sig$, the results 
explained in the last section imply $\bD_\Phi([\sig]) \cong
D^b(\bP_\Phi([\sig]))\cong D^b(A\mof)$, where
$A$ is the endomorphism ring of an injective generator of
$\cP_\Phi(X_\Delta)$, and the grading on
$A$ is given by lifting an injective generator of $\cP_\Phi(X_{[\sig]})$ to 
$\bP_\Phi([\sig])$, taking the endomorphism ring, and using 
Theorem \ref{PC theorem}(3).

Thus to prove Theorem \ref{Koszul functor} and hence Theorem \ref{first main}
it is enough to construct a functor 
\[\kappa\colon \bD_\Phi([\sig])^{opp} \to \bD_{\Phi^\vee}([\sig^\vee])\]
satisfying the properties of Theorem \ref{Koszul functor}.  The grading
shift on $A$-modules becomes the twist $\la k\ra$, and simples and
injective modules in $A\mof$ are replaced by simple and injective
objects of $\bP_\Phi([\sig])$.  Note that we only 
fixed these simples up to a 
twist $\la k\ra$, but property (a) of Theorem \ref{Koszul functor}
implies that there will be only one choice of lift for which 
$\kappa(L^\udot_\tau)$ is perverse.

\subsection{Duality between costandard objects}
To construct the functor $\kappa$ and to prove Theorem 
\ref{third main} we look first at
costandard objects, which have simpler topology than
either injective or simple objects.  If $\tau\prec \sig$, the 
costandard object   
$N^\udot_\tau\in \bP_\Phi([\sig])$ is 
the (derived) ``lower star'' extension of 
the rank one constant local system
$\R_{O_\tau}$ to $X_{[\sig]}$.  
It is a well-known theme in representation theory
that costandards lie midway between simples and
injectives, so it is not surprising that $\kappa$ will send 
costandards to costandards.  In fact, requiring this dictates
our definition of $\kappa$, since costandards generate the 
category $\bD_\Phi(\Delta)$.  

In this way we get a functor 
which is clearly an equivalence of categories, which
has the shift properties of part (a) of Theorem \ref{Koszul functor},
and which satisfies 
$\kappa(N^\udot_\tau) \cong N^\udot_{\tau^\bot}$ for all $\tau\in[\sig]$.  
It remains to show that $\kappa$ interchanges simples and injectives.  
We do this by showing that
the usual ``extend-truncate-repeat'' procedure 
for constructing the simple object 
$L^\udot_\tau$,  becomes under 
$\kappa$ a construction of $I^\udot_{\tau^\bot}$.  
Roughly, the construction of
$I^\udot_{\tau^\bot}$ starts with 
$N^\udot_{\tau^\bot}$ and extends by $N^\udot_\rho$ for
successively smaller $\rho$ (larger strata!), killing off
any possible $\Ext^1$s that remain at each stage.

This follows from the purity of the stalks of $L^\udot_\tau$, which
in our combinatorial language says that the stalk cohomology as a 
bigraded vector space lies all on a diagonal.  We deduce purity 
from the Hard Lefschetz theorem for toric varieties,
which was recently proved for nonrational fans by 
Karu \cite{Ka}, using the theory of equivariant combinatorial sheaves 
on fans of \cite{BrL,BBFK}.  To translate this result to
our category $\bD(\Delta)$ we use the classical Koszul duality
between symmetric and exterior algebras, which relates 
equivariant and ordinary cohomology and sheaf theory \cite{GKM}.

\subsection{Outline of the paper}
We give a brief overview of the sections of this paper.  Section 2
collects a few preliminary definitions regarding gradings,
cones and fans, and sheaves on fans.  Section 3
is concerned with the definition and basic properties of the
categories $\bD(\Delta)$ and $\bD_\Phi(\Delta)$.  The 
categories are defined in \S\S\ref{defining bD} and 
\ref{combinatorial completions}.
The pushforward and pullback functors are described in 
\S\ref{functorial properties}, and the $t$-structure on 
these categories is defined in \S\ref{t-structure}.  The 
injective and costandard objects in the core of this $t$-structure
are studied in \S\S\ref{starshriek} and \ref{costandards and injectives}.
In \S\ref{examples} we work out in detail the structure of these
categories in the simplest nontrivial case of a single one-dimensional
cone.

Section 4 defines the Koszul functor $\kappa$ (\S
\ref{defn of Koszul functor}),
proves the purity of the stalks and costalks of the intersection 
cohomology objects (\S\ref{stalk purity}), and proves 
that $\kappa$ interchanges simples and injectives
(\S\ref{injectives and simples}).

In Section 5, we define the functor $\real_\Delta$, 
and show how the combinatorial functors, completions, $t$-structures,
etc., correspond to topological properties of sheaves.  In particular,
this section contains proofs of Theorems \ref{realization theorem}
and \ref{PC theorem}.

\subsection{Acknowledgments}
This paper has had a long gestation, and there
are many people who should be thanked.  At various points conversations
with R. MacPherson, A. Beilinson, V. Ginzburg, W. Soergel, 
R. Bezrukavnikov, G. Barthel, J.-P. Brasselet, K.-H. Fieseler,
L. Kaup, R. Stanley, K. Vilonen, I. Mirkovic, and V. Lunts have 
been invaluable.

\section{Preliminaries}

\subsection{Gradings, shifts, and total Hom}  \label{gradings}
We begin by fixing notation
regarding gradings.

Let $A = \bigoplus_i A_i$ be a finitely generated
graded $\kk$-algebra, or more generally a sheaf of such algebras on
a finite topological space.
Let  $A\Mof$ and $A\mof$ denote the categories
of finitely generated $A$-modules and finitely generated graded $A$-modules,
respectively.  Let $\{ 1\}\colon A\mof \to A\mof$ be the shift of grading
given by $(M\{ 1\})_i = M_{i+1}$.
Let $\{ n \}$ be the $n$-fold composition of $\{ 1 \}$, or the 
$-n$-fold composition of its inverse if $n < 0$.  

%

For objects $M, N \in A\mof$, we put $\Hom_i(M,N) = 
\Hom_{A\mof}(M,N\{ i\})$, and set 
\[\hom(M,N) = \bigoplus_{i\in \Z} \Hom_i(M,N).\]
If we assume that $\dim_\R A < \infty$, then 
$\hom(M,N)$ is finite dimensional, so it can
be considered as an object in $\R\mof$ with 
$\Hom_i(M,N)$ placed in degree $i$.  

Suppose now that $\cD$ is a full triangulated 
subcategory of either the homotopy category
$K(A\mof)$ or its derived category $D(A\mof)$, and that
$\cD\{ j \} = \cD$ for any $j\in \Z$.  
If $M^\udot$ is a complex in $\cD$, we say that 
elements of $M^i_j$ are in {\em bidegree} $(i,j)$.

Given complexes $M^\udot, N^\udot \in \cD$,
we define
\[\Hom^i_j(M^\udot, N^\udot) = 
\Hom_\cD(M^\udot, N^\udot[i]\{ j \}),\]
\[\hom^\udot(M^\udot, N^\udot) = 
\bigoplus_{i,j\in \Z} \Hom^i_j(M^\udot, N^\udot).\]
We consider $\hom^\udot(M^\udot, N^\udot)$ as an object in 
$D^b(\R\mof)$ with trivial differential.

We have a tensor product functor $D^b(\R\mof)\times \cD \to \cD$
given by 
\[V\otimes M^\udot = \bigoplus_{i,j\in \Z} H^i(V)_j \otimes M^\udot[-i]\{ -j\}.\]
For any $M^\udot, N^\udot\in \cD$, there is a natural morphism
\[\hom^\udot(M^\udot, N^\udot)\otimes M^\udot \to N^\udot.\]

For formal reasons it will be useful to consider complexes $(M^\udot, d)$,
where both the complex and the degree are allowed to be half-integers.  
We let $\Xi$ be the lattice $\{(i,j)\in \R^2 \mid i+j,\, i-j \in \Z\}$, and
we let $K_h(A\mof)$, $D_h(A\mof)$ be the homotopy category and derived category of
the category of complexes $M^\udot = \bigoplus_{(i,j)\in \Xi} M^i_j$ with
differential $d\colon M^i_j \to M^{i+1}_j$ and an action of $A$ so that
$A_kM^i_j \subset M^i_{j+k}$.  All the definitions made above go over to this
setting in an obvious way.

\subsection{Cones} \label{cones}
We fix some notations concerning cones and fans.  A {\em cone} in a 
finite-dimensional real vector space $V$ is a 
subset of the form $\sigma = \R_{\ge 0}v_1 + \dots + \R_{\ge 0}v_k,\;
v_1, \dots, v_k\in V$
which contains no line through $0$.
We include the possibility that $\sig$ is the zero cone $\{0\}$.
The dimension $\dim \sig$ of a cone $\sigma$ is the dimension of 
the vector space $V_\sig = \vspan(\sig)$.  

If a dual vector $\xi \in V^*$ satisfies $\sig \subset \xi^{-1}(\R_{\ge 0})$, 
we call the set $\tau = \sig \cap \xi^{-1}(0)$ a {\em face} of $\sig$.  
It is again a cone in $V$, and we write $\tau \prec \sig$
to indicate the relation that $\tau$ is a face of $\sig$.  If
$[\sig]$ denotes the set of all faces of $\sigma$, 
then $([\sigma],\prec)$ is 
a finite ranked lattice, with maximal element $\sigma$ and minimal element
$o = \{0\}$. 

\begin{rmk}
Note that several of our constructions, including the definition of 
toric varieties, involve taking the dual cone 
$\sig^\vee$ to $\sig$ or the annihilator $V^\bot_\sig\subset V^*$ of the space 
$V_\sig$.  These notions obviously depend on the ambient space
$V$, and we consider the choice of ambient vector space to be a part
of the definition of a cone.  
\end{rmk}
 
\subsection{Fans}
A {\em fan} $\Delta$ in $V$ is a finite collection of cones in 
$V$ for which
\begin{enumerate}
\item[(1)] $\sig \cap \tau$ is a face of both $\sig$ and $\tau$ 
whenever $\sig$ and $\tau$ are in $\Delta$, and
\item[(2)] for each $\sig \in \Delta$, we have $[\sig]\subset \Delta$.
\end{enumerate}
A {\em quasifan} is a finite collection of cones $\Delta$
satisfying (1) and
\begin{enumerate}
\item[($2'$)] if $\sig, \tau \in \Delta$, then 
$\tau \prec \rho \prec \sig$ implies $\rho\in\Delta$.
\end{enumerate}
Clearly every fan is a quasifan.
A quasifan $\Delta$ has a unique smallest fan containing it, namely
$[\Delta] = \bigcup_{\sig\in\Delta} [\sig]$.

Given a quasifan $\Delta$, we give it the topology generated by 
basic open sets $[\sig]\cap \Delta$, $\sig\in \Delta$. 
Quasifans $\Sig\subset \Delta$ are just the locally closed subsets.


Suppose that $V$ contains a lattice $V_\Z$ for which 
$V = V_\Z \otimes_\Z\R$.  A subspace $W\subset V$ is called 
{\em rational} if $W \cap V_\Z$ spans $W$.  A 
quasifan $\Delta$ is said to be rational if
$V_\sig$ is rational for every $\sig\in\Delta$; we call a 
cone $\sig$ rational if $[\sig]$ is rational.


\subsection{Sheaves on fans} \label{sheaf formalities}
Let $\Delta$ be a quasifan in $V$, with the topology described above.
\begin{prop} \label{Sheaves via stalks} 
There is an equivalence of categories
between the category of sheaves of abelian groups on $\Delta$
 and the category of data $(\{S_\sig\},\{r_{\tau,\sig}\})$,
where \begin{itemize}
\item $S_\sig$ is an abelian group for all $\sig \in \Delta$,
\item For every pair of faces $\tau\prec\sig$ in $\Delta$, 
$r_{\tau,\sig}\colon S_\sig \to S_\tau$ is a homomorphism, and
\item whenever $\rho\prec\tau\prec\sig$, we have $r_{\rho,\tau}r_{\tau,\sig}
= r_{\rho,\sig}$.
\end{itemize}
A morphism $\phi\colon (\{S_\sig\}, \{r_{\tau,\sig}\}) \to 
(\{S'_\sig\}, \{r'_{\tau,\sig}\})$ in this category 
is a collection of maps
$\phi_\sig\colon S_\sig \to S'_\sig$ for every $\sig \in \Delta$ 
satisfying $\phi_\tau r_{\tau,\sig} = r'_{\tau,\sig}\phi_\sig$
for every $\tau\prec\sig$.
\end{prop}
\begin{proof}
Given a sheaf $S$, let 
$S_\sig = \Gamma(S;[\sig]\cap \Delta)$; this is the stalk of $S$ at $\sig$,
since $[\sig]\cap \Delta$ is the smallest open set containing $\sig$. 
The homomorphisms $r_{\tau,\sig}$ are then given 
by restriction of sections $\Gamma(S;[\sig]\cap\Delta) \to 
\Gamma(S;[\tau]\cap \Delta)$.

The inverse functor is also easy to describe: 
given data $(\{S_\sig\},\{r_{\tau,\sig}\})$
and an open set $U \subset \Delta$, the space of sections $\Gamma(S;U)$ is 
the inverse limit 
$\lim\limits_{\stackrel{\longleftarrow}{\sig \in U}} S_\sig$. 
\end{proof}
We will pass freely between these points of view on sheaves without comment.

The same statement holds for sheaves of rings, replacing 
``abelian group'' with ``ring'' everywhere.  
Furthermore, given a sheaf $\cR$
of rings on $\Sigma$, we can describe the category
$\cR\mof$ of sheaves
of $\cR$-modules as the category of pairs $(\{S_\sig\},\{r_{\tau,\sig}\})$,
where each $S_\sig$ is an $\cR_\sig$-module and 
each $r_{\tau,\sig}\colon S_\sig\to S_\tau$ is a homomorphism
of $\cR_\sig$-modules.  Here $S_\tau$ becomes an 
$\cR_\sig$-module via $\cR_\sig \to \cR_\tau$.

We also need to describe pullback and pushforward functors in this
language.  Let $\Sig$ be a sub-quasifan of $\Delta$, with inclusion
$j$. If $\cR = \cR_\Delta$ is a sheaf of rings on $\Delta$, the restriction
$\cR_\Sig = \cR_\Delta|_\Sig$ has the same stalks and restriction maps,
but only for $\sig, \tau\in \Sig$.  
In the same way, given an object 
$S \in \cR_\Delta\mof$, the pull-back $j^*S$ is given by keeping only
the stalks $S_\sig$ for $\sig\in\Sig$.

If $S\in \cR_\Sig\mof$, the pushforward sheaf $S' = j_*S$ is given 
by $S'_\sig = \Gamma(S;[\sig]\cap \Sig)$ for $\sig \in \Delta$. 
This is an $\cR_\sig$-module, since it is the inverse limit
of $S_\tau$ for $\tau\in[\sig]\cap \Sig$, each of which
is an $\cR_\sig$-module by the homomorphism $\cR_\sig \to \cR_\tau$.

Finally, we define an ``extension by zero'' functor $j_{!!}$ in this
language (we use this nonstandard notation to avoid confusion with
its derived version, defined in \S\ref{functorial properties}).
For $S \in \cR_\Sig\mof$,  we let $(j_{!!}S)_\tau = 
S_\tau$ if $\tau\in \Sig$, and  $(j_{!!}S)_\tau = 0$ otherwise.

\section{Combinatorial sheaves on fans}

\subsection{The basic definition}\label{defining bD}
Given a quasifan $\Delta$ in $V$, define a sheaf of graded rings 
$\Lambda = \Lambda_\Delta$
on it as follows.  Let the stalk $\Lambda_\sig$ be the exterior algebra
$\Lambda(V_\sig^\bot)$, with the standard grading where 
$\Lambda^1(V_\sig^\bot) = V_\sig^\bot$ has degree $1$. 
If $\tau \prec \sig$ are cones in $\Delta$, let the restriction
$r_{\tau,\sig}\colon\Lambda_\sig \to \Lambda_\tau$ be the natural 
homomorphism induced by the inclusion $V_\sig^\bot \to V_\tau^\bot$.

Let $\LDmod$ denote the abelian category of finitely generated
$\Lambda_\Delta$-modules.  
We write $A\mapsto A\{ k \}$ for the automorphism of  $\LDmod$ 
which shifts the grading down by $k\in \Z$.

For a cone $\sig \in \Sigma$, let $i_\sig\colon \{\sig\} \to \Sig$
denote the inclusion, and define $J_\sig = i_{\sig *}\Lambda_\sig$.
Using Proposition \ref{Sheaves via stalks}, $J_\sig$ has a simple 
description: if $\sig \prec \tau$, then $(J_\sig)_\tau= \Lambda_\sig$
with the induced action of $\Lambda_\tau$, while
$(J_\sig)_\tau=0$ if $\sig \not\prec\tau$.  
The restriction map $r_{\tau,\rho}$ is the identity
whenever $\sig \prec \tau \prec \rho$.
\begin{prop}
The category $\LDmod$ has enough injectives.  An object in 
$\LDmod$ is injective if and only if it is injective as a sheaf
of $\R$-modules.  $J_\sig$ is injective for every $\sig\in\Delta$, 
and any injective object is a direct sum of $J_\sig\{ k\}$
for various $\sig \in\Delta$, $k\in \Z$.
\end{prop}
\begin{proof}
The injectivity of $J_\sig$ follows from the adjunction between 
$i^*_\sig$ and $i_{\sig*}$.  For any $S\in \LDmod$, 
the natural homomorphism
\begin{equation}\label{enough injectives}
 S \to \bigoplus_{\sig\in \Delta} i_{\sig *}i^*_\sig S
\end{equation} 
is injective, so $S$ embeds into the injective object
$\bigoplus_{\sig\in\Delta} i_{\sig *} I_\sig$, where for each
$\sig$ the $\Lambda_\sig$-module $I_\sig$ is an injective hull
of $i_\sig^*S$.
\end{proof}

Let $Inj_\Delta$ be the full subcategory of injectives in $\LDmod$.
We define a triangulated category $\bD(\Delta)$ 
to be the homotopy category $K^b_h(Inj_\Delta)$
of bounded complexes of objects in $Inj_\Delta$,
with the ``half grading'' described in \S\ref{gradings}.  This is
equivalent to the full subcategory of $D^b_h(\LDmod)$ of objects whose
cohomology sheaves are locally free, where $S\in \LDmod$
is called locally free if each stalk $S_\sig$ is a free $\Lambda_\sig$-module. 

Define a functor $\Gamma_\sig\colon\bD(\sig)\to D^b_h(\R\mof)$ 
by 
\[\Gamma_\sig(S^\udot)^p_q = \text{$q$th graded piece of } 
             H^p(S^\udot \otimes_{\Lambda_\sig} \R),\]
with trivial differential.  In terms of our 
dictionary with sheaves on the toric variety $X_\Delta$, 
$\Gamma_\sig(S^\udot)$ should be thought of as the stalk of 
$S^\udot$ at a point of $O_\sig$.

\subsection{Combinatorial completions and stable sheaves}
\label{combinatorial completions}
Let $\Delta$ be a fan in the vector space $V$.  Recall from 
\S\ref{intro to pcs} that a {\em combinatorial completion}
$\Phi$ of $\Delta$ is a choice for each $\tau \in \Delta$
of a complement $\Phi_\tau$ to $V^\bot_\tau$ in $V^*$ satisfying
$\Phi_\rho \subset \Phi_\tau$ whenever $\rho \prec \tau$.
Such completions always exist.  We can, for instance, 
choose a nondegenerate bilinear form on $V^*$
and take $\Phi_\tau$ to be the orthogonal complement of 
$V^\bot_\tau$.

Fix a completion $\Phi$ of $\Delta$.
Given faces $\tau\prec\sig$ of $\Delta$, we define
$\Phi^\tau_\sig = V^\bot_\tau \cap \Phi_\sig$.  It is 
a complement to $V^\bot_\sig$ in $V^\bot_\tau$.
Easy linear algebra shows that
for any cones $\rho\prec\tau\prec\sig$, we have
\begin{equation} \label{Phi direct sum}
\Phi_\sig^\rho = \Phi_\sig^\tau \oplus \Phi_\tau^\rho.
\end{equation}

\begin{defn}
The category $\LDmod_\Phi$ of $\Phi$-{\em stable} sheaves has as
objects all pairs $(S, \{S^\st_\sig\}_{\sig\in \Delta})$, 
where $S \in \LDmod$ and the {\em stabilization}
$(S^\st_\sig)$ is a choice for each $\sig \in \Delta$ of a 
graded vector subspace of $S^\st_\sig \subset S_\sig$ for which:
\begin{enumerate}
\item For each $\sig \in \Delta$, 
multiplication gives an isomorphism 
\[\Lambda_\sig \otimes S^\st_\sig \to S_\sig\] 
(in particular, $S_\sig$ is a free
$\Lambda_\sig$-module), and
\item if $\tau\prec\sig$ are faces in $\Delta$, then 
\[r_{\tau,\sig}(S^\st_\sig) \subset \Lambda(\Phi^\tau_\sig) \cdot 
S^\st_\tau.\]
\end{enumerate}
A morphism $\phi\colon (R,\{R^\st_\sig\}) \to (S,\{S^\st_\sig\})$ 
is a morphism $\phi\colon R\to S$
for which $\phi(R^\st_\sig) \subset S^\st_\sig$ for all $\sig \in \Delta$.
\end{defn}

Forgetting the stabilization gives a functor $\LDmod_\Phi \to \LDmod$. 

\begin{prop} The category $\LDmod_\Phi$ is abelian, and
has enough injectives.
The indecomposable injective $J_\tau$ in $\LDmod$ has 
a unique stabilization 
given by $(J_\tau)^\st_\sig = \Lambda(\Phi^\tau_\sig)$ if
$\tau \prec \sig$, and $0$ otherwise.  The resulting objects are,
up to isomorphism and shifts of grading, the only indecomposable injectives.
\end{prop}

Most of the time we will simply write $S$ instead of $(S,\{S_\sig^\st\})$ 
when discussing objects of $\LDmod_\Phi$.  In particular we will use the
same symbol $J_\tau$ to denote an injective object in $\LDmod_\Phi$
and in $\LDmod$.

Let $Inj_{\Delta,\Phi}$ be the full subcategory of injectives in 
$\LDmod_\Phi$, and define the triangulated category $\bD_\Phi(\Delta)$ to
be the homotopy category $K^b_h(Inj_{\Delta,\Phi})$.  Forgetting the 
stabilizations produces a functor $F_\Delta\colon 
\bD_\Phi(\Delta)\to \bD(\Delta)$.

Intuitively, by passing to the category $\bD_\Phi(\Delta)$ we have made the 
``strata'' $\sig\in\Delta$ contractible.  
The following lemma makes this
precise.

\begin{lemma} \label{contractible strata}
The composed functor 
\[\bD_\Phi(\sig) \stackrel{F_\sig}
{\longrightarrow} \bD(\sig) \stackrel{\Gamma_\sig}
{\longrightarrow} D^b_h(\R\mof)\] is an equivalence of categories.
\end{lemma} 

A more precise description of the structure of the category 
$Inj_{\Delta,\Phi}$ will be useful later.  Let 
$\hom_{\Delta,\Phi} = \hom_{Inj_{\Delta,\Phi}}$ in what follows.

\begin{prop} \label{Homs between injectives} For any faces 
$\sig,\tau\in \Delta$
there is a canonical identification 
\begin{equation}\label{Hom(Jsig,Jtau)}
\hom_{\Delta,\Phi}(J_\tau,J_\sig) = 
\left\{\begin{matrix} \Lambda(\Phi_\sig^\tau)^* & 
\text{if $\tau\prec\sig$,}\\
                     0 & \text{ otherwise}
\end{matrix}
\right.\end{equation}
as graded vector spaces.
If $ \rho\prec \tau\prec \sig$, then 
under the identification {\rm (\ref{Hom(Jsig,Jtau)})}, the composition map
\[\hom_{\Delta,\Phi}(J_\tau,J_\sig) \otimes 
\hom_{\Delta,\Phi}(J_\rho,J_\tau) \to
\hom_{\Delta,\Phi}(J_\rho,J_\sig)\]
is the inverse transpose of the 
wedge product
$\Lambda(\Phi^\tau_\sig) \otimes \Lambda(\Phi^\rho_\tau) \to
\Lambda(\Phi^\rho_\sig)$.
\end{prop}
Note that 
(\ref{Phi direct sum}) shows the wedge product above is an isomorphism.
In (\ref{Hom(Jsig,Jtau)}) the grading is 
given by $\Lambda(\Phi_\sig^\tau)^* = 
\hom_{\R\mof}(\Lambda(\Phi_\sig^\tau),\R)$. In other words, we let elements of
$(\Phi_\sig^\tau)^*$ live in degree $-1$.

\begin{proof}  The identification (\ref{Hom(Jsig,Jtau)}) is given by
\[
\begin{split}
\hom_{\Delta,\Phi}(J_\tau,J_\sig) \to 
\hom_{\sig,\Phi}(i_\sig^*J_\tau, i_\sig^*J_\sig) & \to 
\hom_{\R\mof}((J_\tau)_\sig^\st, (J_\sig)_\sig^\st) \\ & = 
\hom_{\R\mof}(\Lambda(\Phi_\sig^\tau),\R).
\end{split}\]
Checking that this is an isomorphism and 
identifying the composition homomorphism is easy.
\end{proof}

\subsection{Pushforward and pullback functors}   \label{functorial properties}
Fix a quasifan $\Delta$ and a sub-quasifan $\Sigma$; let 
$i\colon \Sigma\to \Delta$ be the inclusion.  Let $\Phi$ be a completion 
of $\Delta$; we will use the same symbol $\Phi$ to denote its restriction to 
$\Sigma$.  In this section, we 
define functors
\begin{equation}\label{eight functors}
\xymatrix{
\bD(\Sig) \ar@<.5ex>[r]^{i_*,i_!} & \bD(\Delta),\ar@<.5ex>[l]^{i^*,i^!}
& \bD_\Phi(\Sig) \ar@<.5ex>[r]^{i_*,i_!} & 
\bD_\Phi(\Delta)\ar@<.5ex>[l]^{i^*,i^!}
}\end{equation}
satisfying familiar properties from sheaf theory.

We begin by remarking that any locally free object $S \in \LDmod$ has a 
bounded functorial resolution
$S \to I^\udot(S)$ by injective objects; just iterate
the construction of (\ref{enough injectives}).

We define the first four functors in (\ref{eight functors}) as follows:
\begin{itemize}
\item The extension functor $i_*$ sends injectives to injectives, so we define
$i_*$ of a complex $J^\udot$ to be the complex $i_*J^\udot$.  
\item 
The restriction functor 
$|_\Sig\colon \LDmod \to \LSmod$ does not take injectives to 
injectives, so it must be derived: let $i^*$ of a complex $S^\udot$ be 
the total complex of the double complex $I^\udot(J^\udot|_\Sig)$.
(Note that if $\Sig$ has a unique smallest cone, it is not 
necessary to derive this functor; this is the only case we will 
actually need.)
\item Similarly, we define $i_!$ by deriving the 
extension by zero functor $i_{!!}\colon \LSmod \to \LDmod$.
In other words, we set $i_!J^\udot = $ 
total complex of $I^\udot(i_{!!}J^\udot)$.
\item
To define $i^!$, let $\ol{\Sigma}$ be the closure of $\Sigma$ in $\Delta$.  
Given a sheaf $S \in \LDmod$ and an open subset $U\subset \Delta$
define $\Gamma_{\ol{\Sig}}(S;U)$ to be the sections of $S$ on $U$ with support 
contained in $\ol{\Sig}$.  The map $U\mapsto \Gamma_{\ol\Sig}(S;U)$ defines 
a sheaf $S_{\ol\Sig}$ on $\Delta$, and the functor $S\mapsto S_{\ol\Sig}$ 
sends injectives to injectives.
Then for a complex $J^\udot$ in $\bD(\Sig)$, let  
$i^!J^\udot = i^*J^\udot_{\ol\Sig}$.
\end{itemize}

Suppose $\Delta$ has a completion $\Phi$ and $\Sigma$ has the
induced completion.  It is easy to check that the functors
$i_*$, $i_{!!}$, $|_\Sig$, $(\cdot)_{\ol\Sig}$, and $I^\udot$ all 
extend in an obvious way to $\Phi$-stable sheaves and complexes, 
giving the other four functors in (\ref{eight functors}).

The functors $i^*$, $i^!$, $i_!$, $i_*$ and their stable versions satisfy the 
following standard properties for pushforward and pullback 
functors (see \cite{GM}); we leave their proof as an exercise for the reader:
\begin{enumerate}
\item[(1)] All four functors are exact functors of triangulated categories.
\item[(2)] If $\Pi\subset \Sigma$ with inclusion $j$,
there are natural isomorphisms $j^*i^* \simeq (i\circ j)^*$, 
$i_*j_* \simeq (i\circ j)_*$, and similarly for 
$i^!$, $i_!$.
\item[(3)] $i^*$ is left adjoint to $i_*$, and $i^!$ is right adjoint to 
$i_!$.  
\item[(4)] There is a natural transformation
$i^! \to i^*$.  If $\Sigma \subset \Delta$ is open, it is an isomorphism.
\item[(5)] There is a natural transformation
$i_! \to i_*$.  If $\Sigma \subset \Delta$ is closed, it is an isomorphism.
\end{enumerate}
For the remaining properties, assume $\Sigma$ is closed in $\Delta$, let
$\Pi = \Delta \setminus \Sigma$, and let $j\colon 
\Pi \to \Delta$ be the inclusion. 
\begin{enumerate}
\item[(6)] We have $j^!i_* = j^*i_* = 0$.
\item[(7)] There are functorial distinguished triangles for $S^\udot$ in 
$\bD(\Delta)$ or $\bD_\Phi(\Delta)$:
\[j_!j^*S^\udot \to S^\udot \to i_*i^*S^\udot \to j_!j^*S^\udot[1]\]
\[i_*i^!S^\udot \to S^\udot \to j_*j^*S^\udot \to i_*i^!S^\udot[1]\]
\item[(8)] The adjunction morphisms
\[i^*i_*S^\udot \to S^\udot \to i^!i_*S^\udot\]
\[j^*j_*\wt{S}^\udot \to \wt{S}^\udot \to j^*j_!\wt{S}^\udot\]
are isomorphisms if $S^\udot\in \bD(\Sig)$ or $\bD_\Phi(\Sig)$, 
$\wt{S}^\udot \in \bD(\Pi)$ or $\bD_\Phi(\Pi)$.
\end{enumerate}

\subsection{The local $t$-structure} \label{t-structure}
To simplify notation, we put $\bD(\sig) = \bD(\{\sig\})$ and
$\bD_\Phi(\sig) = \bD_\Phi(\{\sig\})$ for any cone $\sig$.

The properties (1)--(8) of the previous section 
are exactly what is needed to construct
a $t$-structure 
on $\bD(\Delta)$ by gluing ``local'' $t$-structures 
on the categories $\bD(\sig)$ for all $\sig\in \Delta$.
For more details on $t$-structures and gluing, see
\cite{BBD}, \cite{GM}, \cite{KS}.

To define these local $t$-structures, take $\sig\in \Delta$, and
let $\bD^{\ge 0}(\sig)$ (respectively $\bD^{\le 0}(\sig)$) be 
the full subcategory of objects $S^\udot$ in $\bD(\sig)$ for
which $\Gamma_\sig(S^\udot)^p_q = 0$ if 
$p+q < - \codim(\sig)$ (resp. if $p+q > - \codim(\sig)$).
\begin{prop} 
This defines a $t$-structure on $\bD(\sig)$.
\end{prop}
\begin{proof}
Checking most of the axioms for a $t$-structure
is fairly routine.  We describe the hardest part, 
which is to show that for any $S^\udot \in \bD(\sig)$ there 
is a distinguished triangle
\[S^\udot_{\le 0} \to S^\udot \to S^\udot_{\ge 1} \to S^\udot_{\le 0}[1]\]
with $S^\udot_{\le 0} \in \bD^{\le 0}(\sig)$ and 
$S^\udot_{\ge 1} \in \bD^{\ge 1}(\sig) = \bD^{\ge 0}(\sig)[-1]$.
To show this, suppose that $S^\udot = (S^p_q, d)$, 
and let $S^\udot_{\le 0}$ be the
subcomplex generated under the action of $\Lambda_\sig$ by
\begin{enumerate}
\item all elements of $S^p_q$, if $p+q \le - \codim(\sig)$, and 
\item elements of $S^p_q$ which are in the image of $d$, if  
$p+q = - \codim(\sig) + 1$.
\end{enumerate}
We then let $S^\udot_{\ge 1}$ be the quotient complex 
$S^\udot/S^\udot_{\le 0}$.  It is easy to see that $S^\udot_{\le 0}$ is 
a complex of free $\Lambda_\sig$-modules, so $S^\udot_{\ge 1}$ is also.
The exact sequence
\[0 \to S^\udot_{\le 0} \to S^\udot \to S^\udot_{\ge 1} \to 0\]
induces the required distinguished triangle.  
\end{proof}

The core $\bP(\sig) = \bD^{\ge 0}(\sig)\cap \bD^{\le 0}(\sig)$ of this
$t$-structure is an artinian abelian category.  If $k\in \Z$, the
functor $\langle k\rangle = [k/2]\{ -k/2\}$
is a $t$-exact automorphism of $\bD(\sig)$, and hence induces an 
automorphism of $\bP(\sig)$.  All simple objects of $\bP(\sig)$ 
are of the form $\Lambda_\sig[\codim \sig]\la k \ra$, $k\in\Z$.

Given a completion $\Phi$, 
the $t$-structure on $\bD_\Phi(\sig)$ is defined by using the 
same vanishing restrictions for $\Gamma_\sig F_\sig$.  The fact that
this gives a $t$-structure is trivial, using Lemma
\ref{contractible strata}.

\subsection{The global $t$-structure}
Given a quasifan $\Delta$, define a $t$-structure on $\bD(\Delta)$ by
\[S^\udot \in \bD^{\ge 0}(\Delta) \iff
i^!_\sig(S^\udot) \in \bD^{\ge 0}(\sig)
\;\text{for all}\;\sig \in \Delta,\;\text{and}\]
\[S^\udot \in \bD^{\le 0}(\Delta)\iff 
i^*_\sig(S^\udot) \in \bD^{\le 0}(\sig)
\;\text{for all}\;\sig \in \Delta.\]
The $t$-structure on $\bD_\Phi(\Delta)$ is defined similarly.
The fact that these are $t$-structures follows from the standard
arguments in \cite{BBD} and properties (1)--(8) from 
\S\ref{functorial properties}.

The cores of these $t$-structures are abelian categories, which
we denote $\bP(\Delta)$, $\bP_\Phi(\Delta)$.  Note that 
an object  $S^\udot$ is in $\bP_\Phi(\Delta)$ if and only if its 
image under $F_\Delta$ is in $\bP(\Delta)$.
For any $k\in \Z$, $\la k \ra = [k/2]\{ -k/2\}$
is a $t$-exact automorphism of both $\bD(\Delta)$ and 
$\bD_\Phi(\Delta)$, and 
thus restricts to automorphisms of $\bP(\Delta)$ and $\bP_\Phi(\Delta)$.

\begin{prop} \label{full subcategory}
The forgetful functor $F_\Delta$ makes $\bP_\Phi(\Delta)$ into a 
full subcategory of $\bP(\Delta)$.
\end{prop}

Before proving this, we introduce two useful spectral sequences.
Given objects $A_0^\udot,B_0^\udot \in \bD(\Delta)$, the $\Hom$-spaces
$H^k = \Hom^{k}_{\bD(\Delta)}(A_0^\udot,B_0^\udot)$ are the cohomology
groups of a complex $C^\udot$, with $C^k = \bigoplus_{l\in \Z} 
\Hom_{Inj_\Delta}(A_0^l, B_0^{k+l})$.  Filtering this complex by the 
support of the maps gives a spectral sequence with $E_1$ term
\begin{equation}\label{spectral sequence}
E_1^{p,q} = \bigoplus_{\dim \sig = p} \Hom^{p+q}_{\bD(\sig)}
(i_\sig^* A_0^\udot, i_\sig^! B_0^\udot)
\end{equation}
which abuts to $H^{p+q}$.

In a similar way, if $A^\udot, B^\udot \in \bD_\Phi(\Delta)$, we have
a spectral sequence with $E_1$ term
\begin{equation}\label{spectral sequence 2}
{}_\Phi E_1^{p,q} = \bigoplus_{\dim \sig = p} \Hom^{p+q}_{\bD_\Phi(\sig)}
(i_\sig^* A^\udot, i_\sig^! B^\udot)
\end{equation}
which abuts to $H^{p+q}_\Phi = \Hom^{p+q}_{\bD(\Delta)}(A^\udot,B^\udot)$.

If $A^\udot_0 = F_\Delta A^\udot$, $B^\udot_0 = F_\Delta B^\udot$, then
there is a natural map ${}_\Phi E^{p,q}_\udot \to E^{p,q}_\udot$ 
of spectral sequences, whose $E_1$ component is induced by the
forgetful functors $F_\sig$, and which
abuts to the map $F_\Delta\colon H_\Phi^{p+q} \to H^{p+q}$.

\begin{proof}[Proof of Proposition \ref{full subcategory}] 
Suppose $A^\udot, B^\udot \in \bP_\Phi(\Delta)$, and
take $A_0^\udot$, $B_0^\udot$ as above.  Since 
$i_\sig^*A^\udot \in \bD_\Phi^{\le 0}(\sig)$ and $i_\sig^!B^\udot 
\in \bD_\Phi^{\ge 0}(\sig)$, we have  $E_1^{p,q} = 0$ and 
${}_\Phi E_1^{p,q} = 0$ if $p+q <0$.  As a result, it is enough to 
prove that if $\sig\in \Delta$, 
$A^\udot \in \bD_\Phi^{\le 0}(\sig)$, and $B^\udot \in \bD^{\ge 0}_\Phi(\sig)$,
then the natural map $F_\sig\colon
\Hom^k_{\bD_\Phi(\sig)}(A^\udot, B^\udot) \to 
\Hom^k_{\bD(\sig)}(F_\sig A^\udot, F_\sig B^\udot)$ 
is an injection for all $k$ and an isomorphism if $k = 0$.  

The injectivity follows immediately from
Lemma \ref{contractible strata}.   For surjectivity when $k = 0$, 
note that by replacing $A^\udot$, $B^\udot$ by quasi-isomorphic complexes,
we can assume that $A^d$ is generated in degrees $\le -\dim(\sig) - d$, 
while $B^d$ is generated in degrees $\ge -\dim(\sig) - d$.  
\end{proof}

\subsection{The middle extension} \label{starshriek}
By standard arguments of \cite{BBD}, we have
perverse extension functors
\[j_{!*}\colon \bP(\Sig) \to \bP(\Delta),\;
j_{!*}\colon \bP_\Phi(\Sig) \to \bP_\Phi(\Delta)\]
defined for any inclusion $j\colon \Sigma \to \Delta$ of quasifans.
They satisfy $(ij)_{!*} = i_{!*}j_{!*}$ for any inclusion 
$i\colon \Delta\to \Pi$ of quasifans.

Perverse extension takes simple objects to simple objects.  For 
$\sig \in \Delta$, define 
$L^\udot_\sig = i_{\sig !*}\Lambda_\sig[c/2]\{c/2\}$, where $c =\codim(\sig)$;
since $\bP_\Phi(\Delta)$ is a full subcategory of $\bP(\Delta)$, 
we use the same notation in both categories.  
All simple objects are
isomorphic to $L^\udot_\sig\la k \ra$, for some 
$\sig \in \Delta$ and $k\in \Z$.

There are two descriptions of $j_{!*}$ which will be useful.
If $\Sig$ is closed in $\Delta$, then $j_{!*} = j_! = j_*$, so it is enough to 
let $\Sig$ be an open set.  Let $\Pi$ be the complementary closed set, and
let $i\colon \Pi \to \Delta$ be the inclusion.
\begin{prop} \label{middle extension}
$j_{!*}$ is uniquely characterized by the following properties:
\begin{enumerate}
\item $j^*j_{!*}$ is the identity functor. 
\item For any $S^\udot\in \bP(\Sig)$,  $i^!j_{!*}S^\udot
\in \bD^{\le -1}(\Pi)$ and 
 $i^*j_{!*}S^\udot\in \bD^{\ge 1}(\Pi)$.
\end{enumerate}
\end{prop}

To give a more constructive description of $j_{!*}$, we restrict 
further to the case where $\Pi = \{\sig\}$ is a single cone. 
Then $j_{!*}S^\udot[1]$ is the cone of the 
composed morphism
\[ j_*S^\udot \to i_{\sig *}i^*_\sig j_*S^\udot \to  
i_{\sig *}\tau_{\ge 0}i^*_\sig j_* S^\udot.
\]

\begin{rmk}
The object $L^\udot_o$ constructed in this way 
is essentially
a graded version of the combinatorial intersection cohomology 
complex constructed by McConnell in
\cite{M}, although he expresses it as a sheaf of complexes rather
than a complex of sheaves.  His construction also works for more general 
perversities.
\end{rmk}

\subsection{Costandard and injective perverse objects}
\label{costandards and injectives}
Next we discuss costandard and injective objects in $\bP_\Phi(\Delta)$.
The costandard objects are particularly interesting, since although they
contain many simple perverse objects in their composition series,
when considered as complexes of injectives, they are as simple as possible.

Take $\sig\in \Delta$, and let $c = \codim \sig$.  We define
$N^\udot_\sig = j_{\sig *}\Lambda_\sig[c/2]\{c/2\}$.  
It is simply the injective object
$J_\sig \in Inj_{\Delta,\Phi}$, placed in degree $-c/2$
and grading $-c/2$.  It is easy 
to see that it is perverse, since $j^!_\tau N^\udot_\sig = 0$ if
$\sig\ne \tau$, while if $\sig \prec \tau$, 
then $j_\tau^*N^\udot = \Lambda_\sig[c/2]\{c/2\}$, and   
$\Lambda_\sig$, considered as a $\Lambda_\tau$-module, 
is generated in degrees $0 \le d \le \dim \tau - \dim \sig$.
We call the twists $N^\udot_\sig\la k\ra$, $k\in \Z$, 
{\em costandard} objects.

Dually, we also have {\em standard} objects $M_\sig^\udot\la k\ra$, where
$M^\udot_\sig = j_{\sig !}\Lambda_\sig[c]$.  They are also perverse;
one way to see this would be to define a Verdier duality functor
on $\bD_\Phi(\Delta)$ and show that it preserves perversity and 
interchanges $j_!$ and $j_*$.  It is easy to check perversity directly, 
however, since a minimal complex representing $M^\udot_\sig$ has
\[M^i_\sig \cong \bigoplus_{\stackrel{\sig\prec\tau}{\codim \tau = -i}}
j_{\tau_*}\Lambda_\sig,\]
which implies the required vanishing condition on $j^!_\tau M^\udot_\sig$, while
$j^*_\tau M^\udot = 0$ if $\tau \ne \sig$.

\begin{thm}\label{enough perverse injectives}
 The category $\bP_\Phi(\Delta)$ has enough
injectives; 
each injective object has a filtration by costandard objects.
\end{thm}

This result is essentially the same as \cite[Theorem 3.2.1]{BGS}.
In \S\ref{injectives and simples} we will give an algorithm
based on \cite{BGS} which constructs the indecomposable injective objects
as successive extensions of costandards.  
We do not give the proof that this algorithm works, as it is 
identical to the proof from \cite{BGS}, with two minor modifications.
First, their proof is stated for projective
objects and standard objects rather than injective objects and costandards,
so one must dualize everywhere.  Second,
their proof works in an ungraded setting in which
there is only one standard and one costandard object for each element 
in a finite poset, whereas our standards and costandards can each be
twisted by $\la k\ra$ for any $k\in \Z$; instead of extending by a 
single costandard, therefore, one must extend by all possible 
twists.  

All but one of the hypotheses
of \cite[Theorem 3.2.1]{BGS} hold for formal reasons in any category of
perverse sheaves in which lower star and lower shriek extensions
from single strata are $t$-exact.  
The only one which requires explanation is the following.
\begin{lemma} For any $\tau,\rho \in \Delta$ and $k\in \Z$, we have
\[\Ext^2_{\bP_\Phi(\Delta)}(M^\udot_\tau, N^\udot_\rho\la k\ra) = 0.\]
\end{lemma}
\begin{proof} First note that there is an injection 
(see \cite[Lemma 3.2.4]{BGS})
\[\Ext^2_{\bP_\Phi(\Delta)}(M^\udot_\tau, N^\udot_\rho\la k\ra) \to 
\Hom^2_{\bD_\Phi(\Delta)}(M^\udot_\tau, N^\udot_\rho\la k\ra),\]
so it will be enough to show that this second group is zero.  This 
is equivalent to showing 
$\Hom^2(i_\rho^*M^\udot_\tau, i_\rho^*N^\udot_\rho\la j \ra)=0$.
If $\tau \ne \rho$, we have $i_\rho^*M^\udot_\tau = 0$, so 
we are done. If $\tau = \rho$, use Lemma
\ref{contractible strata}. 
\end{proof}

\begin{thm}
There is an equivalence of categories 
\[\bD_\Phi(\Delta) \cong D^b(\bP_\Phi(\Delta)),\]
in particular, 
$\Hom^j_{\bD_\Phi(\Delta)}(A, B) = \Ext^j_{\bP_\Phi(\Delta)}(A, B)$ 
for any $A, B \in \bP_\Phi(\Delta)$.
\end{thm}
\begin{proof} See \cite[Corollary 3.3.2]{BGS}.
\end{proof} 

It is also true 
that $\bD(\Delta) \cong D^b(\bP(\Delta))$, but the proof is more
involved, and we don't need it here.

\subsection{Example} \label{examples}
We will illustrate the perverse $t$-structure on 
$\bD(\Delta)$ and $\bD_\Phi(\Delta)$ in the simplest 
nontrivial case.
Let $V = \R$, and let $\sig = \R_{\ge 0}$,
$\Delta = [\sig] = \{o,\sig\}$.  This is the fan of the 
toric variety $X_\Delta = \C$, with the orbit stratification
$O_\sig = \{0\}$, $O_o = \C^*$.  


The reader may want to compare the statements that follow with the 
note \cite{BGSch} by Beilinson, Ginzburg, and Schechtman, which
explains the Koszul duality of \cite{BGS} in the special case
$\mathfrak{ g = sl}_2$, where the flag variety is the projective
line $\PP^1$, with stratification $\{\{0\}, \PP^1 \setminus \{0\}\}$.

There are two indecomposable injectives in $\LDmod$, 
up to shifts of grading:
$J_o = i_{o *}\Lambda_o$ and $J_\sig = i_{\sig *}\Lambda_\sig$.
The automorphism groups of both these objects are 
the scalars $\R$, and we have 
 \[\Hom(J_o,J_o\{1\}) \cong 
\Hom(J_o, J_\sig) \cong \Hom(J_o, J_\sig\{-1\})\cong \R.\]
For any other $\rho, \tau \in \Delta$ and $k\in \Z$,
we have $\Hom(J_\rho,J_\tau\{k\}) = 0$.

\begin{table}[tbh] 
\[\begin{array}{|c|l|c|l|} \hline
\;\;\;\;\;\;\;S^\udot & S^{-1} \to S^0  & \real_\Delta S^\udot & 
\xymatrix{
A_o \ar@<.4ex>[r]^{p} & A_\sig\ar@<.4ex>[l]^{q}
}
\\  \hline
L^\udot_\sig = N^\udot_\sig & 
0 \to J_\sig & i_{\sig!*}\R_\sig & \;\xymatrix{
0 \ar@<.4ex>[r]^{0} & \R\ar@<.4ex>[l]^{0}
}\\ 
N^\udot_o\la 1\ra = I^\udot_o\la 1\ra & 
J_o \to 0 & i_{o*}\R_o[1] & \xymatrix{
\R \ar@<.4ex>[r]^{0} & \R\ar@<.4ex>[l]^{1}
}\\
L^\udot_o\la 1\ra & J_o \to J_\sig\{-1\} & i_{o!*}\R_o[1] & \xymatrix{
\R \ar@<.4ex>[r]^{0} & 0\ar@<.4ex>[l]^{0}
}\\
M^\udot_o\la 1\ra & J_o \to J_\sig \oplus J_\sig\{-1\} & 
i_{o!}\R_o[1]\;  & \xymatrix{
\R \ar@<.4ex>[r]^{1} & \R\ar@<.4ex>[l]^{0}
}\\ 
 I^\udot_\sig & J_o \to J_\sig & \cI & \xymatrix{
\R \ar@<.4ex>[r]^{\binom{0}{1}} & \R^2\ar@<.4ex>[l]^{(1,0)}
}\\
\hline
\end{array}\]
\caption{Some perverse objects}
\end{table}

In the first two columns of Table 1 we list five objects in 
$\bP(\Delta)$, along with the names we have been using for them
(the use of $I^\udot_\sig$ and $I^\udot_o$ is potentially misleading,
as these objects are not injective in 
$\bP(\Delta)$; they are the images of the corresponding injectives
in $\bP_\Phi(\Delta)$ under $F_\Delta$).
In the second column, the left element in each complex is 
placed in degree $-1$, and all maps 
have maximal possible rank. 

In the third column we list the corresponding perverse sheaves on 
$X_\Delta$ obtained by applying the functor $\real_\Delta$
of Theorem \ref{realization theorem}.
Here $\R_\tau = \R_{O_\tau}$, the constant rank one local system
on $O_\tau$, for $\tau = o,\sig$. We abuse 
notation and write $i_\tau$ for the inclusion of $O_\tau$
into $X_\Delta$.  

The perverse sheaf
$\cI$ is the largest indecomposable 
extension of $\R_{o}[1]$; it fits in an exact sequence
\[0 \to i_{\sig!*}\R_\sig \to \cI \to i_{o*}\R_o[1] \to 0.\]
This and other standard exact sequences involving these perverse sheaves
can be easily constructed in $\bD(\Delta)$, as the reader is invited to 
check.

In the fourth column we describe these perverse sheaves in
terms of the following well-known result.
\begin{thm}[\cite{MV, V}] \label{perverse sheaves on a line}  
The category of perverse sheaves on $\C$ constructible with 
respect to the stratification $\{\{0\},\C^*\}$ is equivalent to the
category of finite-dimensional vector spaces and maps 
\[\xymatrix{
M_o \ar@<.5ex>[r]^{p} & M_\sig\ar@<.5ex>[l]^{q}}
\]
for which $1+qp$ is invertible.  Under this equivalence, the
vector space $M_o$ is the stalk cohomology at a point of $O_o = \C^*$,
and the operator $1 + qp$ gives the monodromy action of a loop 
around $0$.
\end{thm}

There is a nonzero nilpotent map $\cI\to \cI$.  
It is the image under $\real_\Delta$ of the following map of complexes
(all maps have maximal possible rank):

\begin{equation*}
\xymatrix{I_\sig^\udot\ar[d] & 0 \ar[r]\ar[d] & J_o \ar[r]\ar[d] & 
J_\sig\ar[d] \\
 I_\sig^\udot\la 2\ra &  J_o\{-1\} \ar[r] & J_\sig\{-1\} \ar[r] & 0}
\end{equation*}
It is easy to see that this is not chain-homotopic to $0$.

\subsection{Example of stable sheaves} 
The fan $\Delta$ from the previous section 
has a unique completion given by $\Phi_o = \{0\}$, $\Phi_\sig = V^*$.
Passing from $\bD(\Delta)$ to $\bD_\Phi(\Delta)$ 
corresponds intuitively to embedding $X_\Delta = \C$ into 
$\PP^1$ with strata $O'_\tau = \ol{O_\tau}$, $\tau = o,\sig$. 
These new strata have the same local properties as before, 
but now the open orbit is contractible.

An orbit-constructible perverse sheaf $P$ on $\C$ extends 
to a perverse sheaf on $\C\PP^1$ constructible with respect 
to this stratification 
if and only if $H^{-1}(P|_{O_o})$ is a trivial local
system.  There are five isomorphism classes of indecomposable
perverse sheaves on $X_\Delta$ with this property, namely the five
in Table 1.  On the other hand, the objects in $\bP(\Delta)$
from the left column of the table all have a unique lift to 
$\bP_\Phi(\Delta)$, and they are the only indecomposable objects 
which lift.

The category $\bP_\Phi(\Delta)$ is a mixed version of the category
$\cO$ for $\mathfrak{g = sl}_2$.  The statement of Koszul duality 
in this case is identical with \cite{BGS}.

Note that the nonzero map $I^\udot_\sig \to I^\udot_\sig\la 2\ra$
described above lifts uniquely to a nonzero map in $\bD_\Phi(\Delta)$. 
Since $\sig$ is isomorphic to its own dual cone $\sig^\vee$, the
properties of the Koszul duality functor asserted by Theorem
\ref{Koszul functor} imply that there should be a nonzero map
\[L^\udot_o \to L^\udot_o[2]\la -2 \ra = L^\udot_o[1]\{1\}\]
in $\bD_\Phi(\Delta)$.
Indeed, we have the following map of complexes:
\[\xymatrix{ 
L^\udot_o\ar[d] & 0 \ar[r]\ar[d] & J_o \ar[r]\ar[d] & J_\sig\{-1\} \ar[d] \\
L^\udot_o[1]\{1\} & J_o\{1\} \ar[r] & J_\sig \ar[r] & 0
}\]
Note that this map becomes zero in $\bD(\Delta)$, so 
$F_\Delta\colon \bD_\Phi(\Delta) \to \bD(\Delta)$ is not a full embedding,
although its restriction to $\bP_\Phi(\Delta)$ is. 

\section{Koszul duality}
\subsection{The dual cone} Fix a cone $\sigma\subset V$.
The dual cone $\sigma^\vee\subset 
V^*$ is given by
\[\{\xi \in V^*\mid \xi(v) \ge 0 \text{ for all } v\in \sigma\}.\]
There is an order-reversing bijection $\tau \mapsto \tau^\bot$
between faces of $\sigma$ and faces of $\sigma^\vee$, defined by
\[\tau^\bot = \sigma^\vee \cap (V_\tau)^\bot.\]


Fix a completion $\Phi$ on the fan $[\sig]$, as described
in \S\ref{combinatorial completions}.  It induces a dual completion
$\Phi^\vee$ on $[\sig^\vee]$, by $\Phi^\vee_{\tau^\bot} = (\Phi_\tau)^\bot.$
If $\sig$ is a rational cone, then
 $\Phi$ is rational if and only if $\Phi^\vee$ is rational.

\subsection{The Koszul functor}\label{defn of Koszul functor} Stated in terms of
the combinatorial categories of \S3, our main theorem is the following.

\newcommand{\Omv}{{\Phi^\vee}}
\newcommand{\Sigv}{{\Sig^\vee}}
\newcommand{\sigv}{{\sig^\vee}}
\newcommand{\taub}{{\tau^\bot}}

\begin{thm}  \label{third main}
Given an $n$-dimensional pointed cone $\sig \subset V \cong \R^n$ 
and a combinatorial completion 
$\Phi$ of $[\sig]$, there is an equivalence of 
triangulated categories
\[\kappa\colon \bD_\Phi([\sig])^{opp} \to \bD_{\Phi^\vee}([\sig^\vee])\]
for which 
\begin{enumerate}
\item[(a)]  
$\kappa(M^\udot[k]) = (\kappa M^\udot)[-k]$ and 
$\kappa(M^\udot\langle k\rangle) = 
(\kappa M^\udot)[-k]\langle k \rangle$ for any 
$M^\udot \in \bD_\Phi([\sig])$, and 
\item[(b)] 
$\kappa(L^\udot_\tau) \cong I^\udot_{\tau^\bot}$ and
$\kappa(I^\udot_\tau) \cong L^\udot_{\tau^\bot}$ for any $\tau\prec \sig$.
\end{enumerate}
Here $I^\udot_\tau$, $I^\udot_{\tau^\bot}$ are injective hulls of
the simple objects $L^\udot_\tau$ and $L^\udot_{\tau^\bot}$ in $\bP_\Phi([\sig])$ and
$\bP_{\Phi^\vee}([\sig^\vee])$, respectively.
\end{thm}

We will define this functor as 
the derived functor of a contravariant functor
\[K= K_\sig\colon \Lsmod_\Phi \to 
\Lambda_{[\sig^\vee]}{\mof}_{\Phi^\vee}.\]

Let $k$ denote the functor $\kk\mof \to \kk\mof$
which reverses the grading, so $k(M)_j = M_{-j}$.  
Note that for any vector space $V$, interior 
multiplication makes $k(\Lambda(V^*))$ into a graded $\Lambda(V)$-module. 

Define the functor $K$ as follows.  Take $S\in \Lsmod_\Phi$.
The stalk of $KS$ at a cone $\tau^\bot \in [\sig^\vee]$ is 
given by 
\[(KS)_{\tau^\bot} = k(\Lambda(\Phi_\tau) \otimes \hom(J_\tau, S)^*),\]
where $\hom = \hom_{\Lsmod_\Phi}$ is the graded $\hom$ of \S\ref{gradings}
and $\otimes$ denotes the tensor product of graded  
$\R$-modules.
It has an action of $\Lambda_{\tau^\bot} = \Lambda(V_\tau)$
by interior multiplication on the first 
factor and the trivial action on the second.
It acquires a stabilization by setting 
\[(KS)_{\tau^\bot}^\st =  
k(\omega_\tau \otimes \hom(J_\tau, S)^*),\]
where $\omega_\tau = \Lambda^{\dim \tau}(\Phi_\tau)\cong \R\{-\dim \tau\}$.

To make these stalks into a sheaf, we must define restriction maps
\begin{equation*}\label{sss} 
r_{\rho^\bot,\tau^\bot}\colon (KS)_{\tau^\bot} \to (KS)_{\rho^\bot}
\end{equation*}
for $\tau \prec \rho$. Note that 
Proposition \ref{Homs between injectives} says that 
$\hom(J_\tau,J_\rho) = \Lambda(\Phi^\tau_\rho)^*$ canonically, so 
we can dualize the composition homomorphism to get a map
\begin{equation}\label{ttt}
\hom(J_\tau,S)^* \to \Lambda(\Phi^\tau_\rho) \otimes 
\hom(J_\rho,S)^*.
\end{equation}
Exterior multiplication gives a map (in fact an isomorphism)
\begin{equation}\label{uuu}
\Lambda(\Phi_\tau) \otimes \Lambda(\Phi^\tau_\rho)  
\to \Lambda(\Phi_\rho).
\end{equation}
The required map $r_{\rho^\bot,\tau^\bot}$ is obtained by 
combining (\ref{ttt}) and (\ref{uuu}).  It is a homomorphism
of $\Lambda_{\tau^\bot}$-modules, since $\Lambda_{\tau^\bot}$
acts trivially on $\Lambda(\Phi^\tau_\rho)$.

\begin{rmk}
One can also construct $K$ directly
using the description of morphisms in $Inj_{[\sig],\Phi}$ given by 
Proposition \ref{Homs between injectives}.  
\end{rmk}

\begin{thm} \label{functor K} $K$ is a contravariant functor 
$\Lsmod_\Phi \to \Lsvmod_{\Phi^\vee}$.  
It satisfies the following properties:
\begin{enumerate}
\item $K$ is left exact; we have $K(S\{n\}) = (KS)\{n\}$ for all sheaves $S$ 
and all $k \in \Z$.  
\item For any $\tau\in[\sig]$, there is an isomorphism
\[ KJ_\tau\cong J_{\tau^\bot}\{\dim \tau\}.\]
\item $K$ restricts to an isomorphism of categories
\[Inj_{[\sig],\Phi} \to Inj_{[\sig^\vee],\Phi^\vee}.\]
\end{enumerate} 
\end{thm} 

\begin{proof} 
Using Proposition \ref{Homs between injectives}, it is easy to see
that $r_{\rho^\bot,\tau^\bot} r_{\tau^\bot,\mu^\bot} = r_{\rho^\bot,\mu^\bot}$
whenever $\mu\prec\rho\prec\tau$, so $KS$ is a sheaf.  
It is stable, since
\[r_{\rho^\bot,\tau^\bot}((KS)_{\tau^\bot}^\st) \subset 
k((\omega_\tau \cdot \Lambda(\Phi^\tau_\rho)) \otimes \hom(J_\rho, S)^*) = 
k(\Lambda((\Phi^\vee)^{\rho^\bot}_{\tau^\bot})\cdot 
(KS)_{\rho^\bot}^\st),\]
where the equality follows  
since $(\Phi^\vee)^{\rho^\bot}_{\tau^\bot} =
\Phi^\vee_\tau \cap V_\rho$ is the annihilator of $\Phi_\tau$ in
$V_\rho$.

Part (1) of the theorem is clear.  To show (2), take
$\rho^\bot \in [\sig^\vee]$.  If $\rho\not\prec \tau$, 
then $(KJ_\tau)_{\rho^\bot} = 0$.
If $\rho\prec\tau$, then Proposition \ref{Homs between injectives}
gives a canonical identification
\begin{equation}\label{ppp}
\begin{split}
 (KJ_\tau)_{\rho^\bot} & = 
k(\Lambda(\Phi_\rho) \otimes \hom(J_\rho,J_\tau)^*)\\
& \cong k(\Lambda(\Phi_\rho)\otimes \Lambda(\Phi^\rho_\tau))
\cong k(\Lambda(\Phi_\tau)).
\end{split}
\end{equation}
The direct sum decomposition $V^* = \Phi_\tau \oplus V^\bot_\tau$
shows that $k(\Lambda(\Phi_\tau))$ is isomorphic to $\Lambda_{\tau^\bot}$
as a $\Lambda_{\rho^\bot}$-module.
It is easy to check that the restriction maps between these stalks are 
isomorphisms. 

The essential surjectivity in (3) follows, so we only 
need to show full faithfulness.  It will be enough to show that for 
$\tau\prec \rho$ the map 
\[\Hom_{\Lsmod}(J_\rho,J_\tau) \to \Hom_{\Lsvmod}(K(J_\tau),K(J_\rho)) 
\] 
is injective, since both sides have the 
same dimension.  But applying $\hom(J_\rho,-)$ to a nonzero map 
$J_\rho \to J_\tau$ results in a nonzero map.
\end{proof}

Now, let $\tilde\kappa\colon \bD_\Phi([\sig]) \to \bD_{\Phi^\vee}([\sig^\vee])$
be the derived functor of $K$: if $(I^\udot, d)$ is a complex 
of injectives, then 
$\tilde\kappa(I^\udot, d) = ((KI)^{-\udot}, Kd)$.  The functor $\kappa$
we want is then $\tilde\kappa$ followed by the twist $\la \dim V\ra$. 

It follows from Theorem \ref{functor K} that $\kappa$ is an 
equivalence of categories.  It is also clear that 
part (a) of Theorem \ref{third main} holds.  The proof of part (b)---that
$\kappa^\vee$ interchanges simple and injective objects---will 
occupy the rest of \S4.  The key ingredient is the following result, which 
follows immediately from Theorem \ref{functor K}(2).

\begin{lemma} \label{kappa on costandards}
For any $\tau \in [\sig]$, there is an isomorphism $\kappa N^\udot_\tau\cong N^\udot_{\tau^\bot}$.
\end{lemma}

\begin{rmk}
This isomorphism is not canonical, and the functor 
$\kappa_{\sig^\vee}\colon \bD_{\Phi^\vee}([\sig^\vee])
\to \bD_\Phi([\sig])$ defined by switching the roles of
$\sig$ and $\sig^\vee$ is not canonically the inverse of
$\kappa$.  The correct inverse functor is 
$\kappa^\vee = \omega_\sig \otimes \kappa_{\sig^\vee}$,
where the one-dimensional vector space $\omega_\sig$ is
in bidegree $(0,0)$.  With this definition, there is a canonical 
isomorphism $\kappa^\vee \kappa \simeq \id$.  
\end{rmk}

\begin{rmk} It is possible to view $\kappa$ as a convolution
functor, analogous to the
geometric Fourier transform.  There is a sheaf $\Xi$ on the poset
$[\sig \times \sig^\vee]$ whose stalk on $(\tau, \rho^\bot)$ is
$\Lambda_\sig = \Lambda(V)$ if $\rho \prec \tau$ and is $0$ otherwise.
The support of $\Xi$ can be thought of as a
``combinatorial conormal variety'' to
the fan $[\sig]$.

Exterior and interior multiplication define commuting actions of 
$p_1^*\Lambda_{[\sig]}$ and $p_2^*\Lambda_{[\sig^\vee]}$ on $\Xi$, where 
$p_1$ and $p_2$ are the projections of $[\sig \times \sig^\vee] = 
[\sig]\times [\sig^\vee]$ onto the first and second factors.
Then $\kappa$ can be written
\[S^\udot \mapsto p_{2*}\cRHom_{p_1^*\Lambda_{[\sig]}}(p_1^*S^\udot, \Xi),\]
where $p_{2*}$ and $\cRHom$ must be defined to take the 
completions $\Phi$, $\Phi^\vee$ into account.

It would be interesting to know if our Koszul duality functor can
be described geometrically as a
similar convolution on filtered $D$-modules.
\end{rmk}

\subsection{Purity of the IC sheaves} \label{stalk purity}
We return for the moment to the case
of a general quasifan $\Delta$ in $V$.  We need the following result about
the simple perverse sheaves $L^\udot_\tau$.

\begin{thm} \label{pointwise purity} 
If $\tau \prec \sig$ are cones in $\Delta$, and $j\in \half\Z$, then
$H^j(i_\sig^*L^\udot_\tau)$ and $H^j(i_\sig^!L^\udot_\tau)$ are 
free $\Lambda_\sig$-modules generated in degree $j$.
\end{thm}

To prove this,
we import results from the theory of equivariant
intersection cohomology of fans, developed independently in 
\cite{BBFK}, \cite{BrL}.
The key result is the Hard Lefschetz theorem proved by Karu \cite{Ka}.  
We recall the basic facts of this theory.

Given a fan $\Delta$, let $\cA = \cA_\Delta$ denote the 
sheaf of $\Delta$-piecewise polynomial functions
on $\bigcup_{\sig\in\Delta} \sig$, graded so that
linear functions have degree two.  In the language of 
\S\ref{sheaf formalities},
the stalk $\cA_\sig$ is the symmetric algebra $S(V_\sig^*)$, with the 
the obvious restriction map $\cA_\sig\to \cA_\tau$ 
if $\tau \prec \sig$.

Since all of these rings are quotients of $A = S(V^*)$, any space of
sections of an $\cA$-module will be an $A$-module.  For any
$A$-module $M$, we define $\ol{M} = M \otimes_A \R$.

Let us fix a cone $\tau\in \Delta$. A  
{\em minimal extension sheaf}\/ for $\tau$ on $\Delta$  
is an $\cA$-module $\cL = \cL^\tau$ satisfying:
\begin{enumerate}
\item (Normalization) $\cL_\tau \cong \cA_\tau\{\codim \tau\}$,
and $\cL_\sig = 0$ if $\tau \not\prec \sig$,
\item (Local freeness) $\cL_\sig$ is a free $\cA_\sig$-module for each 
$\sig\in \Delta$, and
\item (Minimal local extension) for each $\sig \in \Delta$, $\sig\ne \tau$, 
the homomorphism
\[\ol{\cL_\sig} \to \ol{\cL(\partial\sig)}\]
is an isomorphism of graded vector spaces.
\end{enumerate}
Here we define $\partial\sig = (\Delta\cap [\sig]) \setminus \{\sig\}$.  

Such a sheaf always
exists, and is unique up to a scalar isomorphism.  The shift $\codim \tau$ 
in (1) does
not appear in \cite{BBFK} or \cite{BrL}; we have added it so the minimal
extension sheaf is perverse.

Properties (1) and (3) imply that $\cL_\sig \to \cL(\partial\sig)$
is surjective for any $\sig$, 
or in other words, that 
forgetting the $\cA$-module structure makes $\cL$ injective as
a sheaf of vector spaces.  Define $\cL(\sig,\partial\sig)$ to be the kernel
of $\cL_\sig \to \cL(\partial\sig)$.
\begin{thm}[\cite{BBFK,Ka}] \label{equivariant HL}
 For every $\sig\in \Delta$, $\sig \ne \tau$,
we have:
\begin{enumerate}
\item[(a)] $\cL(\sig,\partial\sig)$ is a free $\cA_\sig$-module,
\item[(b)] $\ol{\cL_\sig}$ vanishes in degrees $\ge -\codim \sig$, and
\item[(c)] $\ol{\cL(\sig,\partial\sig)}$ vanishes in degrees 
$\le - \codim \sig$.
\end{enumerate}
\end{thm}
\begin{proof} (a) is a special case of Corollary 4.12 of \cite{BBFK}.
(b) is equivalent to (c) by Remark 1.8(ii) of \cite{BBFK}.  Finally,
(b) is equivalent to the ``Hard Lefschetz'' theorem proved by 
Karu \cite{Ka}: the implication HL $\implies$ (b) is given in the
last paragraph of \cite{BBFK}, but their argument runs just as well 
backward.  Note that most of these proofs are stated for the case of
$\tau = o$, but they extend easily to the case of general $\tau$.
\end{proof}

Let $\Lambda^\vee$ be the exterior algebra of $V$, graded so that 
$V$ has degree $-1$.  This ring is the Koszul dual to the 
polynomial ring $A$, in the sense of \cite{BGG,BGS,GKM}.  
We will apply the BGG Koszul functor 
(specifically, the graded version of \cite{BGS}) to the stalks
of $\cL$ to obtain a complex of sheaves of $\Lambda^\vee$-modules.

For a graded $A$-module $M=\bigoplus M_{j}$, 
define a complex $h(M) = (N^\udot, d)$, where 
$N^\udot = M \otimes \Lambda^\vee$
with the grading given by placing $M_{j} \otimes \Lambda^\vee_i$
in complex degree 
$j/2$ and grading degree $j/2-i + \dim V$.  The 
boundary map $d$ is given by choosing a basis $\{v_i\}$ for $V$ and the
dual basis $\{v_i^*\}$ for $V^*$, and setting
\[d(m \otimes \lambda) = \sum_i v_i^*m \otimes (\lambda \wedge v_i).\]

Given our minimal extension sheaf $\cL = \cL^\tau$, 
we define a sheaf of $\Lambda^\vee$-modules
$S^\udot = h(\cL)$ by $S^\udot_\sig = h(\cL_\sig)$. 
The following lemma implies that it has the structure of a complex of 
$\Lambda_\Delta$-modules. 

\begin{lemma}\label{A to Lambda}
Suppose that $M$ is a free $\cA_\sig$-module, so that elements of 
$V_\sig^\bot$ act trivially.  Then 
we can give $h(M)$ the structure of a 
complex of $\Lambda_\sig = \Lambda(V_\sig^\bot)$ 
modules so that  $\Gamma_\sig(h(M))$ 
is isomorphic as a 
bigraded vector space to $\ol{M}$, with $\ol{M}_{j}$ appearing in 
bidegree $(j/2, j/2)$.
\end{lemma}
\begin{proof}
Let an element $a\in \Lambda_\sig$ act via
$a(m \otimes \lambda) = m \otimes i(a)\lambda$,
where $i(a)$ denotes interior multiplication by $a$.
\end{proof}

This gives a $\Lambda_\sig$-module
structure on the stalk $S^\udot_\sig$.  It is clearly compatible with the
maps between faces, so $S^\udot$ becomes a complex of 
$\Lambda_\Delta$-sheaves. 
Since $\cL$ is injective as 
an $\R$-sheaf, so is $S^\udot$.
Using the second part of Lemma \ref{A to Lambda}, we see that although 
$S^\udot$ is only bounded below as a complex of injective 
$\Lambda_\Delta$-sheaves, its cohomology vanishes outside of a
bounded set of degrees, so in fact $S^\udot\in \bD(\Delta)$.  

\begin{lemma} We have $i_\sig^!S^\udot = h(\cL(\sig,\partial\sig))$.
\end{lemma}
\begin{proof} Split $\cL$ into a direct sum of 
indecomposable injective $\R$-sheaves.  The only ones which contribute
to either side are the ones supported on the closure of $\sig$.
\end{proof} 

Combining this with Lemma \ref{A to Lambda} shows that
$S^\udot$ satisfies the same vanishing conditions of 
Proposition \ref{middle extension}, and so 
$S^\udot\cong L^\udot_\tau$.
Theorem \ref{pointwise purity} now follows by using
Theorem \ref{equivariant HL} and Lemma \ref{A to Lambda}.

\begin{rmk} The construction actually gives a functor 
$\cA_\Delta\mof_{pure}\to \bD(\Delta)$, where 
$\cA_\Delta\mof_{pure}$ is the category of {\em pure} $\cA$-modules,
which are finite direct sums of shifts of minimal extension sheaves 
$\cL^\tau$.  In \cite{BL} complexes of pure sheaves are
used to model $T$-equivariant
mixed sheaves on toric varieties in the same way $\bD(\Delta)$ models 
constructible mixed sheaves.
\end{rmk}

\subsection{Injectives and simples} \label{injectives and simples}

Now we prove part (b) of Theorem \ref{third main}, 
which says that $\kappa$ interchanges simple and 
injective objects.  
\begin{blank} {\it Constructing simples.}
We begin by restating the construction of simples from \S\ref{starshriek}
in a suitable form.
Take a quasifan $\Delta$ and 
a face $\tau\in \Delta$ of codimension $c$; we want to construct
the simple perverse object $L^\udot = L^\udot_\tau$. 
Refine the partial order on the cones in the closure 
$\ol{\tau} = \{\rho \in \Delta\mid \tau \prec \rho\}$ to a total 
order; this gives a sequence $\rho_0 = \tau , \rho_1,\dots, \rho_N$,
so that the set 
$\Delta_k = (\Delta \setminus \ol{\tau}) \cup \{\rho_0,\dots,\rho_k\}$
is open in $\Delta$ for $k = 0,\dots, N$. Let 
$j_k\colon \Delta_k\to  \Delta$ be the inclusion.

The inductive construction of $L^\udot$ constructs the objects
$L^\udot_k = j_{k*}j^*_kL^\udot$ for $k = 0,\dots, N$, starting
with $L^\udot_0 = N^\udot_\tau$ and ending with $L^\udot_N= L^\udot$.
Assuming $\Lu_{k-1}$ is already constructed, then
$\Lu_{k}[1]$ is the cone of the map
\[\Lu_{k-1} \to \tau(\hom^\udot(\Lu_{k-1}, 
\Nu_{\rho_k})^*)\otimes \Nu_{\rho_k},\]
where $\tau$ is the truncation on bigraded vector spaces $M$ given by 
$\tau(M)^p_q = M^p_q$ for $p+q \ge 0$ and $\tau(M)^p_q =0$ for $p+q < 0$, and  
the map is obtained from the natural maps
$\Lu_{k-1} \to \hom^\udot(\Lu_{k-1}, \Nu_k)^* \otimes \Nu_k$ and 
$M \to \tau(M)$.
 
So far this is just writing out the inductive computation of
\S\ref{starshriek}.  Using Theorem \ref{pointwise purity} we
can alter the algorithm in the following way.
\begin{thm}\label{weight vs degree truncation}
 The algorithm described above is unchanged if the 
truncation $\tau$ is replaced by $\tau'$, which keeps only those
$M^p_q$ for which $q = p + 1$.
\end{thm} 

\begin{proof}
Let $\rho = \rho_k, r = \codim \rho/2$.  We have 
\[
\begin{split}
\hom^\udot(\Lu_{k-1}, \Nu_k)^* & = 
\hom^\udot(j_{\rho}^*\Lu_{k-1}, \Lambda_\rho\{r\}[r])^*\\
& = \Gamma_{\rho}(j_{\rho}^*\Lu_{k-1})\{-r\}[-r].
\end{split}\]
The theorem thus follows from the following
fact:
$\Gamma_{\rho}(j_{\rho}^*\Lu_{k-1})^p_q = 0$
unless either $q = p < r$ or $q = p + 1 > r$.
To show this, 
let $i_k\colon \Delta\setminus \Delta_{k-1} \to \Delta$ 
be the inclusion.  Applying $j_{\rho}^*$ 
to the triangle
\[i^{}_{k*}i^!_k\Lu \to \Lu \to \Lu_{k-1} \to i^{}_{k*}i^!_k\Lu[1]\]
gives a triangle
\[j^!_{\rho}\Lu \to 
j^*_{\rho}\Lu \to j_{\rho}^*\Lu_{k-1} \to j^!_{\rho}\Lu[1].\]
The first map vanishes
by the vanishing conditions for $\Lu$ 
and Lemma \ref{contractible strata}.  Thus 
\[j_\rho^*\Lu_{k-1} \cong j^*_\rho\Lu \oplus j^!_\rho\Lu[1].\]
The claim now follows from the vanishing conditions for 
$\Lu$ and Theorem \ref{pointwise purity}.
\end{proof} 

\end{blank}
\begin{blank} {\it Constructing injectives.} To complete the proof of
Theorem \ref{third main}(b), we give an algorithm for constructing
injective objects in $\bD_\Phi(\Delta)$ which is dual under the 
functor $\kappa$ to the construction of simples given 
above.  For a proof that the following
algorithm does indeed give an injective object, see \cite[Theorem 3.2.1]{BGS}.

Take a cone $\sig$ and a completion $\Phi$ of
$[\sig]$. Let $\Delta = [\sig]$, and fix a face $\tau\in\Delta$ of dimension $d$. 
Extend the partial order on $[\tau]$ to a total order to get a sequence 
$\tau = \rho_0,\dots,\rho_N = o$ so that $\{\rho_0,\dots,\rho_k\}$
is closed in $[\tau]$ for every $0 \le k \le N$.
We construct a sequence of objects
$\Iu_0,\dots,\Iu_N$, starting with $\Iu_0 = \Nu_\tau$ and
ending with $\Iu_N = \Iu_\tau$, so that $\Iu_k$ is
an injective hull of $\Lu_\tau$ in the Serre subcategory of
$\bP_\Phi(\Delta)$ generated by $\{\Lu_\rho\mid\rho \in 
\{\rho_1,\dots,\rho_k\}\cup (\Delta \setminus [\tau])\}$.

For the inductive step, assume that $\Iu_{k-1}$
has already been constructed.  Set $\rho = \rho_k$.
If $l\in \Z$, let
$E_{l} = \Ext^1_{\bP_\Phi(\Delta)}(\Nu_{\rho}\la l \ra,\Iu_{k-1})$.  The
canonical element in 
$E^*_{l}\otimes E_{l}$ induces an element
of $\Ext^1_{\bP_\Phi(\Delta)}(E_{l} \otimes \Nu_\rho\la l \ra,\Iu_{k-1})$.  
Combining these
over all $l$ gives an element of
\[\Ext^1_{\bP_\Phi(\Delta)}(\bigoplus_{l\in \Z} 
E_{l} \otimes \Nu_\rho\la l\ra, \Iu_{k-1}).\]
Define $\Iu_k$ to be the corresponding extension object.

Since 
$\Ext^1_{\bP_\Phi(\Delta)} = \Hom^1_{\bD_\Phi(\Delta)}$, 
we can restate this as follows:
$\Iu_k$ is the object which fits into the distinguished triangle
\[t(\hom^\udot_{\bD_\Phi(\Delta)}(\Nu_\rho, \Iu_{k-1}))\otimes 
\Nu_\rho \to  \Iu_{k-1} \to \Iu_k \to \dots,\] 
where $t$ is the truncation functor on bigraded vector spaces which keeps
only elements of bidegree $(i, 1 - i)$, $i \in \Z$. 

Comparing with Theorem \ref{weight vs degree truncation}
and using Lemma
\ref{kappa on costandards} gives
Theorem \ref{third main}(b).
\end{blank}

\section{Connection with topology}
In this section, we construct the functor $\real_\Delta
\colon \bD(\Delta) \to D^b(X_\Delta)$ described in the Introduction, 
when $\Delta$ is a rational fan.  We show that the combinatorial 
pushforward functors, stabilizations,
etc., defined in Section 3 correspond under this functor to
standard topological functors.  In particular, we prove
Theorems \ref{realization theorem} and \ref{PC theorem} from 
the Introduction.

Given an object $S^\udot \in \bD(\Delta)$, i.e.\
a complex of sheaves in $Inj_\Delta$, we let
$\real_\Delta(S^\udot)$ be the single complex
associated to the double complex 
$\Omega^\udot_\Delta(S^\udot)$, where  
$\Omega^\udot$ is an additive functor from 
$Inj_\Delta$ to the category $Kom(X_\Delta)$ of 
complexes of $\R$-sheaves on $X_\Delta$.
For this definition to satisfy the 
required properties, we need 
$\Omega^\udot_\Delta(J_\tau)$ to be quasi-isomorphic to 
$i_{\tau *}\R_{O_\tau}$ for each $\tau$, and for
the composition $Inj_\Delta \stackrel{\Omega^\udot_\Delta}{\longrightarrow}
Kom(X_\Delta) \to D^b(X_\Delta)$ to be fully faithful.

There are several ways to produce such a functor.
We use ``controlled differential forms'': forms
which on a neighborhood of a point on a stratum $O_\sig$ are isomorphic to
a pullback by the canonical projection to $O_\sig$.  A similar construction
was used by Goresky, Harder, and MacPherson \cite{GHM} and Saper \cite{Sa}
to construct complexes on compactifications of locally
symmetric spaces.

\subsection{Toric geometry}
We briefly recall the rudiments of toric geometry.  

Let $V$ be a finite-dimensional $\R$-vector space spanned by a 
lattice $V_\Z\subset V$.  Let $V^*_\Z\subset V^*$ be the
dual lattice.  For any $S\subset V^*$, let $S_\Z = S\cap V^*_\Z$.
A rational cone $\sig \subset V$
defines an affine toric variety $X_{[\sig]} = \Spec A_\sig$,
where $A_\sig$ is the ring  
$\bigoplus_{v\in\sig^\vee\cap V^*_\Z} \C\phi_v$ with 
multiplication given by 
$\phi_v\phi_{v'} = \phi_{v + v'}$.

If $\tau \prec \sig$ are rational cones, we get a 
natural inclusion of rings 
$A_\sig \to A_\tau$, which defines
an open embedding $X_{[\tau]} \hookrightarrow X_{[\sig]}$.  
For a rational fan $\Delta$ the associated toric variety 
$X_\Delta$ is the union of 
$X_{[\sig]}$, $\sig\in \Delta$, glued along these inclusions.

The torus $T = \Hom(V_\Z, \C^*)$ acts on $X_\Delta$, with one orbit
$O_\sig$ for each cone $\sig\in \Delta$.  A point $x$ lies in 
$O_\sig$ if and only if 
$\phi_v(x) \ne 0 \iff v \in V_{\sig,\Z}^\bot$.
Restricting from $\sig^\bot_\Z$ to 
$V^\bot_{\sig,\Z}$ gives an isomorphism 
$O_\sig \cong \Hom(V_{\sig,\Z}^\bot, \C^*) 
\cong (\C^*)^{\codim \sig}$. 

If $\tau \prec \sig$, then there is a
quotient map $p_{\tau,\sig}\colon O_\tau \to O_\sig$, given
by 
\[\phi_v(p_{\tau,\sig}(x)) = \left\{\begin{array}{ll}
\phi_v(x) &  \text{if}\; v\in V^\bot_{\sig,\Z} \\ 
0 & \text{if}\; v\in \sig^\vee_\Z \setminus  V^\bot_{\sig,\Z}.
\end{array}\right.\]
If $\rho\prec\tau\prec \sig$, then 
$p_{\tau,\sig}p_{\rho,\tau} = p_{\rho,\sig}$.

If $\Delta$ is a quasifan contained in a fan $\Sig$, then we define
$X_\Delta$ to be $\bigcup_{\sig\in\Delta} O_\sig$, the union of 
$T$-orbits in $X_\Sig$.  This does not depend on the choice of $\Sig$.  

\subsection{Realization functor} \label{real}

Let $\Delta$ be a rational quasifan in a vector space $V$.
For each $\sig \in \Delta$,
the complex $(\Omega^\udot_\sig, d)$ of $\R$-valued $k$-forms on $O_\sig$
is a fine resolution of the constant sheaf $\R_{O_\sig}$.  The torus $T$
acts on global sections of $\Omega^\udot_\sig$ by pullback; we will 
work with $T$-invariant forms $\Gamma(O_\sig, \Omega^\udot_\sig)^T$.

The exponential map ${\exp}\colon V/V_\sig\otimes \C \to O_\sig$,
defined by $\phi_x(\exp(\xi)) = e^{2\pi\sqrt{-1}\la \xi,v\ra}$,
is invariant under translations by the lattice $V_\Z$.  Thus if 
we consider elements of $\Lambda_\sig = \Lambda(V/V_\sig)^*$ to
be constant differential forms on $V/V_\sig$, we get a 
map $\chi_\sig\colon \Lambda_\sig \to
\Gamma(O_\sig, \Omega^\udot_\sig)$ whose image consists of $T$-invariant
forms.

\begin{lemma} \label{formality}
$\chi_\sig$ is injective. It induces a quasi-isomorphism
of dg-algebras $\Lambda_\sig \to \Gamma(O_\sig, \Omega^\udot_\sig)$, 
where 
the differential on $\Lambda_\sig$ is trivial.
If $\tau \prec \sig$ are cones in $\Delta$, then
\[p^*_{\tau,\sig}\circ \chi_\sig = \chi_\tau \circ \iota_{\tau,\sig},\]
where $\iota_{\tau,\sig}\colon \Lambda_\sig \to \Lambda_\tau$ is the
natural homomorphism (i.e., the structure homomorphism for the sheaf
$\Lambda_\Delta$).
\end{lemma} 

Given a free graded $\Lambda_\sig$-module $M$, define a 
complex of sheaves $\Omega_\sig^\udot(M)$ on $O_\sig$ as follows.  
The space of sections of $\Omega_\sig^\udot(M)$ on an open set $U$ is
$\Gamma(U, \Omega^\udot_\sig) \otimes_{\Lambda_\sig} M$,
where the (right) action of $\Lambda_\sig$ on 
$\Gamma(U, \Omega^\udot_\sig)$ is given by 
$\omega\lambda = \omega\wedge \chi_\sig(\lambda)|_U$.  The grading is
given by letting $\omega\otimes m$ have degree $i+j$ when $\omega$
is an $i$-form and $m \in M_j$.  The boundary map is given
by $d(\omega\otimes m) = d\omega \otimes m$.  This is well defined
since the previous lemma shows $d \circ \chi_\sig = 0$.

\begin{lemma} \label{stalks}
If $M \cong \bigoplus_j \Lambda_\sig[n_j]$, then 
$\Omega_\sig^\udot(M) \cong \bigoplus_j\Omega^\udot_\sig[n_j]$.
It is a fine resolution of $\bigoplus_j \R_{O_\sig}[n_j]$.  
\end{lemma}

Let $\LDmod_{\lf} \subset \LDmod$ denote the full subcategory of 
locally free sheaves.  Given 
$S\in\LDmod_{\lf}$, we will define a fine complex of sheaves 
$\Omega_\Delta^\udot(S)$ on $X_\Delta$ by gluing the sheaves
$\Omega^\udot_\sig(S_\sig)$ for $\sig\in \Delta$.
Let us describe the sections of this complex on an open
set $U\subset X_\Delta$.  For any $\sig\in \Delta$, put
$U_\sig = U\cap O_\sig$.

If $\tau\prec\sig$ are cones in $\Delta$, then
we have a map
\[q_{\tau,\sig}\colon
\Gamma(U_\sig, \Omega^\udot_\sig(S_\sig)) 
\to \Gamma(p_{\tau,\sig}^{-1}(U_\sig), \Omega^\udot_\tau(S_\tau)),\]
defined by 
$q_{\tau,\sig}(\omega \otimes s) = 
p_{\tau,\sig}^*\omega \otimes r_{\tau,\sig}(s)$, where 
$r_{\tau,\sig}\colon S_\sig\to S_\tau$ is the structure homomorphism.
Lemma \ref{formality} ensures that this is well defined.
It also clearly commutes with the differential $d$.

Define sections of $\Omega_\Delta^\udot(S)$ on $U$ to be tuples
\begin{equation}\label{Omega sections}
(\alpha_\sig)\in \bigoplus_{\sig\in \Delta} 
\Gamma(U_\sig, \Omega^\udot_\sig(S_\sig))
\end{equation}
which are ``locally pullbacks'' along the projections
$p_{\tau,\sig}$: for
any $\sig \in \Delta$ and any point $x\in U_\sig$, 
there must be a neighborhood $\wt{U}$ of $x$ contained
in $U$, so that for any $\tau \prec \sig$ the sections
$\alpha_\tau$ and $q_{\tau,\sig}(\alpha_\sig)$
agree on $\wt{U}\cap p^{-1}_{\tau,\sig}(U_\sig)$.
We define a boundary map on this complex by $d(\alpha_\sig) = (d\alpha_\sig)$.

\begin{lemma} \label{fine sheaf}
 For any object $S \in \LDmod_{\lf}$, the sheaf $\Omega_\Delta^\udot(S)$
is fine.
\end{lemma}
\begin{proof} Sections of $\Omega_\Delta^\udot(S)$ can be multiplied 
pointwise by sections of $\Omega^0_\Delta(\Lambda_\Delta)$, which are
functions $f\colon X_\Delta \to \R$ which
are locally pullbacks as defined above.
There are partitions of unity of such functions subordinate to any open
cover, so $\Omega_\Delta^\udot(S)$ is fine.
\end{proof}

To define the functor $\real_\Delta$, assume first that
we have $S^\udot \in K^b(Inj_\Delta)$ (with no
``half-grading'').
Since $Inj_\Delta$ is a full subcategory of $\LDmod_{\lf}$,
we can let $\real_\Delta(S^\udot)$ be the single complex 
coming from the double complex 
$\Omega^\udot_\Delta(S^\udot)$.  
A general object in $\bD(\Delta) = K^b_h(Inj_\Delta)$ is of the
form $S^\udot = S_0^\udot \oplus S_1^\udot\la 1\ra$, where each 
$S_i^\udot$ is a complex in $Inj_\Delta$; we let 
$\real_\Delta(S^\udot) = \real_\Delta(S_0^\udot) \oplus 
\real_\Delta(S_1^\udot)$.

Since $\Omega^\udot_\Delta(S\{k\}) = 
\Omega^\udot_\Delta(S)[k]$ for any $S\in \LDmod_{\lf}$ and 
$k\in \Z$, property (1) of 
Theorem \ref{realization theorem} is obvious.  It is also clear
that $\real_\Delta$ takes values in $D^b_u(X_\Delta)$, the full subcategory
of orbit-constructible objects in $D^b(X_\Delta)$ with unipotent monodromy,
since $\real_\Delta(J_\sig)\in D^b_u(X_\Delta)$ for every $\sig\in \Delta$.

\subsection{Functorial properties} Before proving parts (2)--(5) 
of Theorem \ref{realization theorem}, we
prove two intermediate results.  Fix an inclusion of
rational quasifans $j\colon \Sig \to \Delta$.

\begin{lemma} \label{functorial lemma}
\begin{enumerate}

\item For any $S \in \LDmod_{\lf}$, there is a natural isomorphism
\[\Omega^\udot_\Sig(j^*S) \cong j^*\Omega_\Delta^\udot(S).\] 
\item For any $S\in Inj_\Sig$, there is a natural quasi-isomorphism 
\[\Omega^\udot_\Delta(j_*S) \stackrel{\sim}{\to} 
j_*\Omega_\Sig^\udot(S).\]
In particular, $\Omega^\udot_\Delta(J_\sig) 
\cong j_{\sig *}\Omega_\sig^\udot$.
\end{enumerate}   
\end{lemma}

\begin{proof} To see (1), note that the compatibility conditions
among the components of a section $(\ref{Omega sections})$ 
imply that the spaces of germs of sections of $\Omega^\udot_\Delta(S)$ 
and of $\Omega_\sig^\udot(S_\sig)$
at a point $x\in O_\sig$ are isomorphic.

To construct the morphism for (2), consider an open set
$U \subset X_\Delta$.  The total space of sections of 
$j_*\Omega_\Sig^\udot(S)$
is $\Gamma(U\cap X_\Sig, \Omega_\Sig^\udot(S))$.  This is given by
tuples as in (\ref{Omega sections}), satisfying the same 
compatibility condition, but where $\sig$ runs over $\Sig$
instead of $\Delta$.
So if $(\alpha_\sig)_{\sig\in \Delta}$ is a section of 
$\Omega_\Delta^\udot(j_*S)$
on $U$, forgetting $\alpha_\sig$ for $\sig\notin\Sig$ provides the 
required section of $j_*\Omega_\Sig^\udot(S)$.

To show this is a quasi-isomorphism, it is enough to show it induces
an isomorphism on stalk cohomology groups at every point $x\in X_\Delta$.
Furthermore, we can assume that $S = J_\sig$, and in fact  
that $\Sig = \{\sig\}$ and $S = \Lambda_\sig$.  Suppose that $x \in O_\tau$.
Since $\Omega_\Sig^\udot(S)$ is a complex of fine sheaves, 
$j_*\Omega_\Sig^\udot(S)$ is isomorphic to the derived push-forward 
$Rj_*\Omega_\Sig^\udot(S)$ in $D^b(X_\Delta)$.  It follows that
\[H^\udot(j_*\Omega_\Sig^\udot(S))_x \cong \Lambda(V^\bot_\sig/V^\bot_\tau).\]
On the other hand, by (1) and
Lemma \ref{stalks} we
have
\[H^\udot(\Omega_\Delta(j_*S))_x 
\cong \R \otimes_{\Lambda_\tau} S_\tau  = \R \otimes_{\Lambda_\tau}
\Lambda_\sig.\]
The induced map between these stalks is the obvious one induced
by the quotient $V^\bot_\sig \to V^\bot_\sig/V^\bot_\tau$; it is 
an isomorphism.
\end{proof}

\begin{cor} If $J, J' \in Inj_\Delta$, then 
the natural map
\[\Hom_{\LDmod}(J,J') \to \Hom_{D^b(X_\Delta)}(\Omega_\Delta^\udot(J),
\Omega_\Delta^\udot(J'))\]
is an isomorphism (i.e.\ $\Omega_\Delta^\udot$ is fully faithful on $Inj_\Delta$).
\end{cor}
\begin{proof} It is enough to take $J = J_\sig$, $J' = J_\tau\{k\}$
for some $\sig,\tau\in \Delta$ and $k\in \Z$.  Applying $j_\tau^*$
to both sides and using the adjunction with $j_{\tau *}$, we can
reduce to the case where $\Delta = \{\tau\}$ is a single cone and
$J = \Lambda_\tau$, $J' = \Lambda_\tau\{k\}$.  In this case 
the right-hand side becomes 
\[\Hom_{D^b(X_\Delta)}(\Omega^\udot_\tau, \Omega^\udot_\tau[k]) = 
\Hom_{D^b(X_\Delta)}(\R_{O_\tau},\R_{O_\tau}[k]) = H^k(O_\tau).\]
The left-hand side is $\Lambda^k(V^\bot_\tau)$, and the map 
sends $\lambda$ to exterior multiplication by 
$\chi_\tau(\lambda)$, acting on 
$\Omega^\udot_\tau$.  The result now follows from
Lemma \ref{formality}.
\end{proof}

Part (2) of Theorem 
\ref{realization theorem} follows.
The first two isomorphisms of part (3) of Theorem \ref{realization theorem}
(the ones for $j^*$, $j_*$) follow from Lemma 
\ref{functorial lemma}.  

To show that 
$j_! \circ \real_\Sig \simeq \real_\Delta \circ j_!$,
note that Lemma 
\ref{functorial lemma}(1) implies that $\real_\Delta(j_{!!} S^\udot)$
is the extension by $0$ of $\real_\Sig S^\udot$, where $j_{!!}$ is the 
combinatorial extension by $0$ functor from \S\ref{sheaf formalities}.
Since $j_!S^\udot$ and $j_{!!}S^\udot$
are quasi-isomorphic as complexes in $\LDmod_{\lf}$, they remain 
quasi-isomorphic after applying $\real_\Delta$.

Finally we show that there is a functorial isomorphism 
\[ \real_\Sig \circ j^! S^\udot \stackrel{\sim}{\to}
j^! \circ \real_\Delta S^\udot, \] where $S^\udot \in \bD(\Delta)$.
To construct the homomorphism, use the previous paragraph and the 
$j_!$, $j^!$ adjunction, which holds in both $\bD(\Delta)$ and 
$D^b(X_\Delta)$.  It is easy to show it is an isomorphism when 
$S^\udot = J_\sig$; the general case follows by induction.

\subsection{Stalks} 
Recall the functor $\Gamma_\sig\colon \bD(\Delta) \to 
D^b_h(\R\mof)$ defined in \S\ref{defining bD}.  We will show that
it corresponds under $\real_\Delta$ to the stalk cohomology 
functor at any point $x\in O_\sig$.  First, define a ``topological'',
nonmixed version of $\Gamma_\sig$:
\[\Gamma^{top}_\sig(S^\udot)^k = \bigoplus_{p+q = k} 
\Gamma_\sig(S^\udot)^p_q.\]

\begin{thm} Take a cone $\sig$ in the rational fan $\Delta$,
and any point $x\in O_\sig$.  Let $i_x\colon \{x\} \to X_\Delta$ 
denote the inclusion.  Then
there is a natural isomorphism of graded vector spaces
\begin{equation} \label{stalk map}
 H^\udot(i_x^* \real_\Delta S^\udot)\stackrel{\sim}{\to} 
\Gamma^{top}_\sig(S^\udot).
\end{equation}
\end{thm}

This result, together with the results of the previous section, 
immediately implies Theorem \ref{realization theorem}, parts (4) and (5).

\begin{proof}
A cohomology class in $H^k(i_x^* \real_\Delta S^\udot)$ is 
represented by a germ of a section of $\Omega^k_\Delta(S^\udot)$
at $x$ which is killed by the differential $d$.  Such a germ
is represented by a tuple $(\alpha_\sig)$ of the form 
(\ref{Omega sections}).  Suppose 
$\alpha_\tau = \sum_{i = 1}^n \omega_i \otimes s_i$, where each 
$\omega_i$ is a section of $\Omega^\udot_\tau$, and $s_i \in S^\udot_\tau$.  
Send this to $\sum_{i=1}^n \omega_i(x)\ol{s_i}\in 
\Gamma^{top}_\tau(S^\udot)^k$, where  $\omega(x)$
is the value of $\omega$ at $x$ if $\omega$ is a $0$-form, and
is zero otherwise, and $s \mapsto \ol{s}$ is the natural quotient
map $S^\udot_\tau \to \R\otimes_{\Lambda_\tau} S^\udot_\tau$.

It is easy to check that this is a well defined map, 
and that it passes to cohomology, giving a map 
(\ref{stalk map}).  To check that it is an isomorphism it is enough
to handle the case $S^\udot = j_{\sig !}\Lambda_\sig$, which is easy:
the stalk at a point of $O_\tau$ is $\R$ if $\tau = \sig$ and $0$ 
otherwise.
\end{proof}

\subsection{$\Phi$-stable sheaves} \label{Phi-stable}
\newcommand{\hSig}{{\wh{\Sig}}}
 We now prove Theorem \ref{PC theorem}, which 
gives the topological interpretation of 
combinatorial completions and the stable category $\bD_\Phi(\Delta)$.  
Our first goal is to prove the ``only if'' direction of part (1) of 
this theorem: that objects of $\bP_\Phi(\Delta)$ are sent into 
$\cP_\Phi(X_\Delta)$
by $\real_{\Delta,\Phi}$.  To do this, we must study how our 
categories of combinatorial and topological sheaves are affected by 
changing the ambient vector space of a fan.

Let $\Sig$ be a rational quasifan in $V$, and suppose that every cone 
of $\Sig$ is contained in a rational subspace $W \subset V$.
For each $\sig\in \Sig$, we denote by $\hat\sig$ the same cone,
but thought of as a subset of the ambient space $W$ instead of
$V$.  Let $\hSig$ be the fan $\{\hat\sig\mid\sig\in \Delta\}$.
The vector space span of the cones has not changed: we have
$W_{\hat\sig} = V_\sig$.  But
the sheaf of rings $\Lambda_{\hSig}$ 
is different, since $(\Lambda_{\hSig})_{\hat\sig} = 
\Lambda_{\hat\sig} = \Lambda(W_{\hat\sig}^\bot)$, where $W_{\hat\sig}^\bot$
is the annihilator of $W_{\hat\sig}$ in $W$.

Now suppose $Z \subset V^*$ is a rational subspace 
complementary to the annihilator $W^\bot$ of 
$W$ in $V^*$.   The composition $Z \hookrightarrow V^* \to
W^*$ is an isomorphism, so 
we get an inclusion $\iota\colon W^* \to Z \to V^*$.   
For any $\sig\in \Sig$ we have $\iota(\hat\sig^\vee) \subset \sig^\vee$, 
so we get a morphism $\pi_Z\colon X_{\Sig}\to X_{\hSig}$,
defined by $\phi_v(\pi_Z(x)) = \phi_{R(v)}(x)$
(here the toric variety $X_{\hSig}$ is defined with respect to 
the lattice $W^*_{\Z} = \iota^{-1}(V^*_\Z) \subset W^*$).
It is a fiber bundle whose fibers are tori $(\C^*)^{\dim V - \dim W}$.

Next, we define the analogue for combinatorial sheaves of the functor
$\pi^*_Z$.  We have inclusions
$\Lambda_{\hat \sig} \to \Lambda_\sig$, $\sig \in \Sig$, which 
are compatible with the structure maps of the sheaves
$\Lambda_\Sig$ and $\Lambda_{\hSig}$.  Thus we can define a functor 
$G_Z\colon \Lambda_{\hSig}\mof \to \LDmod$ by 
$G_Z(S)_\sig = \Lambda_\sig \otimes_{\Lambda_{\hat\sig}} S_\sig$.
It induces a triangulated 
functor $\bD(\hSig) \to \bD(\Sig)$, which we again denote by 
$G_Z$.

\begin{prop} \label{stability}
There is a natural quasi-isomorphism 
\[\pi_Z^* \circ \real_{\hSig} \simeq \real_{\Sig}\circ\, G_Z\]
of functors $\bD(\hSig) \to D^b(X_\Sig)$.
\end{prop}

\begin{proof} First consider the case where $\Sig = \{\sig\}$ is a 
single cone, so $X_\Sig = O_\sig$, $X_{\hSig} = O_{\hat\sig}$.
Take an object $M \in \LDmod_{\lf}$, i.e.\ a free $\Lambda_{\hat\sig}$-module.
Define a homomorphism 
\begin{equation}\label{qqq} \pi_Z^*\Omega^\udot(M) \to \Omega^\udot(G_Z(M))
\end{equation}
as follows.  A section of $\pi_Z^*\Omega^\udot(M)$ will be of the form 
$\pi^*\omega \otimes m$ for $\omega$ a section of 
$\Omega^\udot_{\hat\sig}$ and $m \in M$.  Send this
to $\pi^*\omega \otimes (1\otimes m)$, which is a section 
of $\Omega^\udot(G_Z(M)) = \Omega^\udot(\Lambda_\sig 
\otimes_{\Lambda_{\hat\sig}}M)$.  The fact that this is well defined
follows from the identity 
\begin{equation} \label{nnn}
\chi_\sig \circ \iota = \pi_Z^*\circ \chi_{\hat\sig},
\end{equation}
 which is essentially Lemma \ref{formality} in another
guise.  To show that it is an isomorphism, it is enough 
to consider $M = \Lambda_{\hat\sig}$, in which case it
follows from (\ref{nnn}) and Lemma \ref{formality}.

For general $\Sig$, the homomorphism (\ref{qqq}) is defined 
componentwise by the same formula.  To show it is an isomorphism,
use the case $\Sig = \{\sig\}$ and Lemma \ref{functorial lemma}.
\end{proof}

Let $\Phi$ be a rational completion 
on a quasifan $\Delta$.
Take a cone $\sig\in\Delta$, and let $\Sig = [\sig] \cap \Delta$.
Let $W = V_\sig$, $Z = \Phi_\sig$.
The previous discussion gives a fan $\hSig$ in $W$ 
and a map $\pi=\pi_\sig\colon X_\Sig\to X_\hSig$; this is the map 
used in Definition \ref{defining stable perverse sheaves} to define 
the category $\cP_\Phi(X_\Delta)$.

\begin{prop} \label{zzz} There is a functor 
$Q\colon \bD_\Phi(\Sig) \to \bD(\hSig)$ 
so that $G_Z\circ Q$  is naturally isomorphic to
the forgetful functor
$\bD_\Phi(\Sig)\to \bD(\Sig)$.
\end{prop}
\begin{proof} Suppose $(S,\{S^\st_\tau\})\in \LDmod_\Phi$.  
If $\hat\tau\in \hSig$, 
define $Q(S)_{\hat\tau} = 
\Lambda_{\hat \tau}S^\st_\sig \subset S_\sig$.  Since
$\Lambda_{\hat\tau} = \Lambda(\Phi^\tau_\sig)$, 
these stalks fit together to give an object in 
$\Lambda_{\hSig}\mof$. 
\end{proof}

\begin{cor} \label{cor}  
 If $S^\udot\in \bD_\Phi(\Delta)$, then 
$\real_{\Delta,\Phi}( S^\udot|_{X_\Sig}) \cong 
\pi^*\real_{\hSig}\wh{S}^\udot$ for some
object $\wh{S}^\udot\in \bD(\hSig)$.
\end{cor} 
\begin{proof} Take $\wh{S}^\udot = Q(S^\udot|_{X_\Delta})$ and use
Propositions \ref{stability} and \ref{zzz}.
\end{proof}

Theorem \ref{PC theorem} (1) follows.  

\subsection{Ext's and injectives}

Finally, we prove the remaining parts of Theorem \ref{PC theorem}. 
To simplify notation, we will indicate the result
of applying $\real_{\Delta,\Phi}$ by a change to 
script letters: $\cL_\sig = \real_{\Delta,\Phi}(L^\udot_\sig)$,
$\cI_\sig = \real_{\Delta,\Phi}(I^\udot_\sig)$, etc.
By the results of the previous section, these objects lie in 
$\cP_\Phi(X_\Delta)$.
 
We begin with a special case of Theorem \ref{PC theorem}(3).  
Define a bifunctor $\ext^1_{\bP_\Phi(\Delta)}$
by \[\ext^1_{\bP_\Phi(\Delta)}(A^\udot,B^\udot) = \bigoplus_{k\in \Z}
\Ext^1_{\bP_\Phi(\Delta)}(A^\udot,B^\udot\la k \ra).\]
\begin{prop} \label{Ext^1 to simple}
Suppose $S^\udot \in \bP_\Phi(\Delta)$ and $\sig\in \Delta$.  Then the map
\begin{equation} \label{total realization}
\ext^1_{\bP_\Phi(\Delta)}(L^\udot_\sig, S^\udot) \to 
\Ext^1_{\cP_\Phi(X_\Delta)}(\cL_\sig, \cS)
\end{equation}
induced by $\real_{\Phi,\Delta}$ 
is an isomorphism.
\end{prop}
Assuming this for the moment, we obtain 
\begin{cor} $\cI_\sig$ is an injective hull of $\cL_\sig$.
\end{cor}
\begin{proof} Taking $S^\udot = I^\udot_\sig$ in Proposition
\ref{Ext^1 to simple} shows that $\Ext^1( \cL_\tau, \cI_\sig) = 0$ for
all $\tau\in \Delta$, so $\cI_\sig$ is injective.  Applying 
$\real_{\Delta,\Phi}$ to the injection
$L^\udot_\sig \to I^\udot_\sig$ gives an injection $\cL_\sig \to \cI_\sig$. 
To see that $\cI_\sig$ is indecomposable, note that 
$\bP_\Phi(\Delta)$ and $\cP_\Phi(X_\Delta)$ are
full subcategories of $\bP(\Delta)$ and $\cP(X_\Delta)$, respectively,
and hence full subcategories of $\bD(\Delta)$ and $D^b(X_\Delta)$.
So using Theorem \ref{realization theorem}(2), we see that if 
$\cI_\sig$ were decomposable, there would be a nontrivial idempotent
in the ring 
\[\bigoplus_{k\in \Z}\Hom_{\bP_\Phi(\Delta)}(I^\udot,
I^\udot\la k \ra).\]  Such an idempotent
must actually lie in $\Hom_{\bP_\Phi(\Delta)}(I^\udot,
I^\udot)$, contradicting the indecomposability of $I^\udot$.
\end{proof}
This implies part (3) of Theorem \ref{PC theorem}: given
$S_1^\udot$, $S_2^\udot \in \bD_\Phi(\Delta)$, take injective 
resolutions and apply Theorem \ref{realization theorem} (2) again.

We can now deduce the ``if'' direction of Theorem \ref{PC theorem} (1).
Suppose $S^\udot\in \bP(\Delta)$, and $\cS = \real_\Delta S^\udot$
is in $\cP_\Phi(\Delta)$.  Since $\cP_\Phi(\Delta)$ has enough 
injectives by the previous corollary, 
there is an injection $\phi\colon \cS \to \cI$, where
$\cI \in \cP_\Phi(\Delta)$ is injective; say 
$\cI = \real_{\Delta,\Phi} I^\udot$, where
 $I^\udot$ is injective in $\bP_\Phi(\Delta)$. 
Replacing $I^\udot$ by
$\bigoplus_{k = -K}^K I^\udot\la k \ra$ if necessary, we can 
construct 
$f \in \Hom_{\bD(\Delta)}(S^\udot, I^\udot)$ which is
sent to $\phi$ by $\real_\Delta$.  Since $\phi$ is an injection, 
so is $f$.  We have embedded $S^\udot$ into an object from
$\bP_\Phi(\Delta)$, so $S^\udot\in\bP_\Phi(\Delta)$ also.

Finally, we prove Proposition \ref{Ext^1 to simple}.
Consider the commutative square
\[\label{square}
\xymatrix{ \ext^1_{\bP_\Phi(\Delta)}(L^\udot_\sig, S^\udot) 
\ar[r]\ar[d]^{R_{\Delta,\Phi}} 
& \ext^1_{\bP(\Delta)}(L^\udot_\sig, S^\udot) \ar[d]^{R_\Delta} \\
\Ext^1_{\cP_\Phi(X_\Delta)}(\cL_\sig, \cS) \ar[r] & \Ext^1_{\cP(\Delta)}
(\cL_\sig,\cS)
}\]
Here $R_{\Delta,\Phi}$ is the map (\ref{total realization}), the other 
vertical map $R_\Delta$ 
is defined similarly, and the horizontal maps are  the injections induced
from the forgetful functors.  We call
an element of $\ext^1_{\bP(\Delta)}(L^\udot_\sig, S^\udot)$ ``$\Phi$-stable''
if it is in the image of the upper horizontal map.

To show 
$R_{\Delta,\Phi}$ is injective, it is enough to show that
$R_\Delta$ is injective.  This follows from Theorem 
\ref{realization theorem}(2), 
since $\Ext^1_{\bP(\Delta)} \to \Hom^1_{\bD(\Delta)}$ and 
$\Ext^1_{\cP(X_\Delta)} \to \Hom^1_{D^b(X_\Delta)}$ are isomorphisms
of bifunctors.

To show that $R_{\Delta,\Phi}$ is surjective, take 
$\psi \in \Ext^1_{\cP_\Phi(X_\Delta)}(\cL_\sig, \cS)$,
with image $\phi \in \Ext^1_{\cP(X_\Delta)}(\cL_\sig, \cS)$.
There exists $\tilde{\phi} \in 
\ext^1_{\bP(\Delta)}(L^\udot_\sig, S^\udot)$ so that 
$R_\Delta(\tilde\phi) = \phi$.  We must show that $\tilde{\phi}$
is $\Phi$-stable. 

To do this, 
use the setup of Proposition \ref{zzz}: let $\Sig = \Delta\cap [\sig]$,
let $\hSig$ be the same quasifan with the ambient space $W = V_\sig$, and 
let $G = G_{\Phi_\sig} \colon \bD(\hSig) \to \bD(\Sig)$.  The shifted 
functor $G' = G[\codim \sig]$ preserves perversity.

\begin{lemma} A map $\tilde\phi \in 
\ext^1_{\bP(\Delta)}(L^\udot_\sig, S^\udot)$ is $\Phi$-stable
if and only if the restriction $\tilde\phi|_{\Sig}$ is in the image 
of $G'\colon \ext^1_{\bP(\hSig)}(L^\udot_\sig, S^\udot)
\to \ext^1_{\bP(\Sig)}(L^\udot_\sig, S^\udot)$.
\end{lemma}
\begin{proof} The ``only if'' direction follows from Corollary \ref{cor}.
For the other direction,  note that 
for any $\tau\in \Delta$ we have 
$i^!_\tau S^\udot \in \bD^{\ge 0}(\Delta)$ and 
$i^*_\tau L^\udot_\sig \in
\bD^{\le a}(\Delta)$, where $a = 0$ if $\tau = \sig$ and $a = -1$
if $\tau\ne \sig$. 
As in the proof of 
Proposition \ref{full subcategory}, the spectral sequences
(\ref{spectral sequence}) and (\ref{spectral sequence 2}) imply  
that if $\tau \ne \sig$, then the map 
\[\hom^k_{\bP_\Phi(\Delta)}(i^*_\tau L_\sig^\udot,i^!_\tau S^\udot) \to 
\hom^k_{\bP(\Delta)}(i^*_\tau L_\sig^\udot,i^!_\tau S^\udot)\]
is an isomorphism for $k = 1$, an injection for $k > 1$, and both sides
vanish if $k \le 0$, while if $\tau = \sig$, it is an isomorphism when
$k=0$, an injection if $k > 0$, and both sides vanish if $k < 0$.

As a result, in order for $\tilde\phi$ to be
$\Phi$-stable,
it is enough for $\tilde\phi|_{\Sig}$ to be $\Phi$-stable.
But $L^\udot_\sig|_{\Sig}$ has support only on $\sig$, where it is
a shifted copy of $\Lambda_\sig$.  So $\Phi$-stability of 
$\tilde\phi|_{\Sig}$ is equivalent to 
$\tilde\phi|_{\Sig}$ being in the image
of $G'$.
\end{proof}

Now, since $\phi$ comes from an extension in $\cP_\Phi(X_\Delta)$, 
$\phi|_{X_\Sig}$ is in the image of 
\[\pi^*\colon \ext^1_{\cP(X_{\hSig})}(\wh\cL_\sig, \wh\cS)
 \to \ext^1_{\cP(X_{\Sig})}(\cL_\sig, \cS)\]
for some $\wh\cS$ and $\wh\cL_\sig \in \cP(X_{\hSig})$.
It follows, using Proposition \ref{stability} and Theorem
\ref{realization theorem}(2) that
$\tilde\phi|_{\Sig}$ is in the image of $G'$, so we can apply
the previous lemma to complete the proof of Proposition
\ref{Ext^1 to simple}.

\end{document}